\documentclass[final]{siamltex}
\usepackage{graphicx,amssymb,amsmath, cancel, lineno}
\usepackage{empheq}

\makeatletter
\def\elsartstyle{%
    \def\normalsize{\@setfontsize\normalsize\@xiipt{14.5}}
    \def\small{\@setfontsize\small\@xipt{13.6}}
    \let\footnotesize=\small
    \def\large{\@setfontsize\large\@xivpt{18}}
    \def\Large{\@setfontsize\Large\@xviipt{22}}
    \skip\@mpfootins = 18\p@ \@plus 2\p@
    \normalsize
}
\@ifundefined{square}{}{}
\makeatother

\MHInternalSyntaxOn
\providecommand*\phantomword[3][c]{%
\mathchoice
{\MT_phantom_word:NNnn #1\displaystyle {#2}{#3}}%
{\MT_phantom_word:NNnn #1\textstyle {#2}{#3}}%
{\MT_phantom_word:NNnn #1\scriptstyle {#2}{#3}}%
{\MT_phantom_word:NNnn #1\scriptscriptstyle {#2}{#3}}%
}
\def\MT_phantom_word:NNnn #1#2#3#4{%
\@begin@tempboxa\hbox{$\m@th#2#4$}%
\setlength\@tempdima{\widthof{$\m@th#2#3$}}%
\hbox{\hb@xt@\@tempdima{\csname bm@#1\endcsname}}%
\@end@tempboxa}
\MHInternalSyntaxOff

\newtheorem{thm}{Theorem}[section]
\newtheorem{Def}{Definition}[section]
\newtheorem{prop}{Proposition}[section]
\newtheorem{Cor}{Corollary}[section]

\newtheorem{Lem}{Lemma}[section]
\newtheorem{rmk}{Remark}[section]

\newcommand{\mc}{\mathcal}

\title{Exact boundary controllability results for a multilayer Rao-Nakra sandwich beam}

\author{A. \"{O}zkan \"{O}zer \thanks{Department of Applied Mathematics, University of Waterloo, Waterloo, ON N2L3G1, Canada ({\tt aozer@uwaterloo.ca}).}
        \and Scott W. Hansen \thanks{Department of Mathematics, Iowa State University,
 Ames, Iowa 50011, USA ({\tt shansen@iastate.edu}).  \quad Supported in part by the National Science Foundation under grant DMS-1312952.}}
\begin{document}

\maketitle

\begin{abstract}
We study the boundary controllability problem for a multilayer
Rao-Nakra sandwich beam. This beam model  consists of a Rayleigh beam coupled
with a number of wave equations.   We consider  all   combinations of clamped and hinged boundary
conditions with the control applied to either the moment or the rotation angle at an end of the beam.  We prove that exact controllability
holds provided  the damping parameter is sufficiently small. In
the undamped case, exact controllability holds without any restriction on the parameters in the system.  In each case, optimal control time is obtained in the space of optimal regularity for $L^2(0,T)$ controls.   A key step in the proof of our main result is the proof of uniqueness of the zero solution of the eigensystem with the homogeneous boundary conditions together with zero boundary observation.
 \end{abstract}



\pagestyle{myheadings}
\thispagestyle{plain}
\markboth{}{Submitted - \"OZER AND HANSEN}

 \enlargethispage*{2.5pc}

 \section{Introduction}
 The classical  sandwich beam is an engineering model for a  three layer beam consisting of two ``face plates" and a ``core" layer that is orders of magnitude more compliant than the face plates.    While most of the early models considered only transverse dynamics, e.g., \cite{Mead-Marcus}, \cite{Yan-Dowell},
 the model due to Rao and Nakra \cite{Rao-Nakra}  includes rotary inertia in each layer and longitudinal inertia (in addition to transverse inertia).
 The model assumes continuous, piecewise linear displacements through the cross-sections, with the Kirchhoff hypothesis imposed on  the face plates.

 In this article we study the boundary controllability of the following  multilayer generalization of the Rao-Nakra beam derived in \cite{Hansen3}:
\begin{equation}\left\{ \begin{array}{l}
m\ddot w -\alpha \ddot w'' +  K   w'''' -   N^{\rm T} {\bf{h}}_E \left({\bf{ G}}_E \psi_E + \tilde{{\bf{ G}}}_E \dot \psi_E \right)' = 0 ~~~ {\rm{in}}~~ \Omega\times \mathbb{R}^+ \\
   {\bf{h}}_{\mc O} {\bf{p}}_{\mc O} {\ddot y}_{\mc O} -{\bf{h}}_{\mc O} {\bf{E}}_{\mc O}  {y}_{\mc O}'' + {\bf{B}}^{\rm T}  \left({\bf{ G}}_E \psi_E + \tilde{{\bf{ G}}}_E \dot \psi_E \right) = 0~~ {\rm{on}} ~~ \Omega\times \mathbb{R}^+ \\
 {\text{where}}~~  {\bf{B}} { y}_{\mc O}={\bf{h}}_E \psi_E-{\bf{h}}_E  N w',
 \end{array} \right.
\label{maincont}
\end{equation}
where $\Omega=(0,L)$, primes denote differentiation with respect to the spatial variable $x$ and dots denote differentiation with respect to time $t$.

The model (\ref {maincont})  consists of  $2m+1$ alternating stiff and complaint (core) layers, with stiff layers on outside.   The stiff layers have odd indices $1,3,\ldots 2m+1$ and the even layers have even indices $2,4,\ldots 2m$.    The Kirchhoff hypothesis is imposed on the stiff layers and Timoshenko displacement assumptions are assumed in the compliant layers.    Damping proportional
rate of shear is included in the compliant layers.

 In the above, $m,\alpha, K$ are \emph{positive} physical constants, $w$ represents the transverse displacement, $\psi^i$ denotes the shear angle in the $i^{\rm{th}}$ layer, $\psi_E=[\psi^2,\psi^4,\ldots,\psi^{2m} ]^{\rm T},$ $y^i$ denote the longitudinal displacement
along the center of the   $i^{\rm{th}}$   layer, and $y_{\mc O}=[y^1,y^3, \ldots, y^{2m+1}]^{\rm T},$   and
\begin{eqnarray}
\nonumber &{\bf{p}}_{\mc O}={\rm{diag}}~ (\rho_1,  \ldots,\rho_{2m+1}),  ~~{\bf{h}}_{\mc O}={\rm{diag}}~(h_1, \ldots, h_{2m+1}), ~~{\bf{h}}_E={\rm{diag}}~(h_2,\ldots, h_{2m}),~~&\\
\nonumber &{\bf{E}}_{\mc O}={\rm{diag}}~(E_1, \ldots, E_{2m+1}), ~~~{\bf{ G}}_E={\rm{diag}}~( G_2, \ldots, G_{2m}), ~~\tilde{{\bf{ G}}}_E={\rm{diag}}~({\tilde{ G}}_2,  \ldots, {\tilde{ G}}_{2m})&
\end{eqnarray}
where $h_i, \rho_i, E_i,$ are positive and denote the thickness, density, and Young's modulus, respectively.
Also $G_i\ge 0$ denotes shear modulus of the $i^{\rm{th}}$ layer, and $ {\tilde{ G}}_i\ge 0$ denotes coefficient for damping in the corresponding compliant layer.

The vector $ N$ is defined as
$N={\bf{h}}_E^{-1}{\bf{A }}{\bf{h}}_{\mc O}  \vec 1_{\mc O} +  \vec 1_E$
where ${\bf{A}}=(a_{ij})$  and ${\bf{B}}=(b_{ij})$ are the $m\times(m+1)$ matrices
$$a_{ij}  = \left\{ \begin{array}{l}
1/2,~~{\rm{  if   }}~~j = i~~{\rm{  or }}~~j = i + 1 \\
~~0,\quad{\rm{   otherwise}} \\
\end{array} \right., ~~b_{ij}=\left\{ \begin{array}{l}
(-1)^{i+j+1},~~{\rm{  if   }}~~j = i~~{\rm{  or }}~~j = i + 1 \\
~~0, \quad\quad\quad\quad {\rm{    otherwise}} \\
\end{array} \right. $$
and $ \vec 1_{\mc O}$ and $\vec 1_E$ denote the vectors with all entries $1$ in $\mathbb{R}^{m+1}$ and $\mathbb{R}^{m},$ respectively.

Consider (\ref{maincont}) with either hinged-Neumann (h-N), or clamped-Dirichlet (c-D), or mixed-mixed (m-m) boundary conditions respectively
\begin{eqnarray}\label{bdryback10}
 &  \left\{ \begin{array}{l}
  w(0,t)=w''(0,t)= w(L,t)=0, w''(L,t) = M(t) ~~ {\rm{on}}~ \mathbb{R}^+\label{bdrycont3}\\
 {y}_{\mc O}'(0,t) = 0, ~ {y}_{\mc O}'(L,t)=  {\bf{g}}_{\mc O}(t)~~~ {\rm{on}}~~ \mathbb{R}^+,\label{bdrycont4}
   \end{array} \right\} &~~~~\text{(h-N)} \\
 &  \left\{ \begin{array}{l}
  w(0,t)=w'(0,t)= w(L,t)=0, w'(L,t) = M(t) ~~~ {\rm{on}}~~ \mathbb{R}^+\label{bdrycont5}\\
 {y}_{\mc O}(0,t) = 0,  ~{y}_{\mc O}(L,t)=  {\bf{g}}_{\mc O}(t) ~~~ {\rm{on}}~~ \mathbb{R}^+,\label{bdrycont6}
   \end{array} \right\} &~~~~\text{(c-D)} \\
 & \left\{ \begin{array}{l}
  w(0,t)=w'(0,t)= w(L,t)=0, w''(L,t) = M(t) ~~~ {\rm{on}}~~ \mathbb{R}^+\label{bdrycont7}\\
 {y}_{\mc O}(0,t) = 0,  ~ {y}_{\mc O}'(L,t)=  {\bf{g}}_{\mc O}(t)~~~ {\rm{on}}~~ \mathbb{R}^+ .\label{bdrycont8}
   \end{array} \right\} &~~~~\text{(m-m)}
\end{eqnarray}
The initial conditions for (\ref{maincont}) are
\begin{eqnarray}
w(x,0)=w^0(x)&,&  ~\dot w(x,0)=w^1(x), ~{y}_{\mc O}(x,0)= {y}^0_{\mc O}, ~ {\dot y}_{\mc O}(x,0)= {y}^1_{\mc O} ~ {\rm{on}}~~ \Omega .\label{initialcont}
\end{eqnarray}

In this paper, through the controls $M(t)$ and ${\bf{g}}_{\mc O}(t)$ at the right end of the  beam, we control the moment and longitudinal force of the stiff layers in (\ref{bdrycont4}) and (\ref{bdrycont8}), and  the shear angle and the longitudinal displacements of the stiff layers in (\ref{bdrycont6}).

\subsection{Background}    In \cite{Rajaram-Hansen2}, exact boundary controllability of three-layer Rao-Nakra beam was investigated for the boundary conditions  (\ref{bdrycont5}).   An exact controllability  result for sufficiently  large control time but with size restrictions on the coupling parameters ($\tilde{ {\bf{ G}} }$ and ${{\bf{ G}}}$ in (\ref{maincont})) was obtained by the standard multiplier method.      In \cite{Rajaram-Hansen3},  the moment method was applied  to the three-layer Rao-Nakra system with  the boundary conditions   (\ref{bdrycont3}).       Under the assumption of distinct wave speeds, exact controllability   was shown up to a finite-dimensional subspace which consists of low-frequency eigenvectors  of the system.   With additional  restrictions on the parameters ($\tilde{ {\bf{ G}} }$ and ${{\bf{ G}}}$ in (\ref{maincont})), and  exact controllability  of the vibrational states was obtained.       Exponential boundary feedback stabilization results for a related (but different) three layer laminated beam were obtained in \cite{Wang}.
In \cite{HO1}, \cite{HO2}  exact controllability results for the multilayer Rao-Nakra plate system analogous to (\ref{maincont}) with locally distributed control in a neighborhood of a portion of the boundary were obtained by the method of Carleman estimates.

\subsection{Main results}

Let \begin{subequations}
\label{control}
\begin{empheq}[left={\phantomword[r]{}{ \mc{C}} = \empheqlbrace}]{align}
\label{control1}& (   H^2(\Omega)\cap H^1_0(\Omega))\times (\tilde H^1(\Omega))^{(m+1)} \times H^1_0(\Omega) \times (\tilde L^2(\Omega))^{(m+1)}   & \text{(h-N)}~\\
\label{control2}&H^1_0(\Omega)\times (L^2(\Omega))^{(m+1)} \times  (L^2(\Omega)/ {\mathrm{M}}) \times (H^{-1}(\Omega))^{(m+1)} & \text{(c-D)}~\\
\label{control3}& H^2_\#(\Omega)\times (H^1_\dag(\Omega))^{(m+1)} \times H^1_0(\Omega)\times (L^2(\Omega))^{(m+1)} & \text{(m-m)}~
\end{empheq}
\end{subequations}
where $\tilde H^1(\Omega)$ and $\tilde L^2(\Omega)$  are the quotient spaces  defined by $\tilde H^1(\Omega)=H^1(\Omega)/\mathbb{R}$ and $\tilde L^2(\Omega)=L^2(\Omega)/ {\mathbb{R}}$ respectively, and
\begin{eqnarray}\nonumber
\nonumber {\mathrm{M}}&&={\rm span}\{e^{-\frac{1}{\sqrt{\alpha/m}}x}, e^{\frac{1}{\sqrt{\alpha/m}}x}\},\\
\nonumber H^2_{\#}(\Omega)&&=\left\{u\in H^2(\Omega)\cap H^1_0(\Omega)~:~ u'(0)=0\right\},\\
\label{spaces} H^1_\dag(\Omega)&&=\left\{u\in H^1(\Omega)~:~ u(0)=0\right\}.
\end{eqnarray}

\begin{prop}\label{weaksolution}Let $T>0, $ and $(M(t),{\bf{g}_{\mc O}}(t))\in (L^2(0,T))^{(m+2)}.$ For any $(w^0,  {y}^0_{\mc O}, w^1, {y}^1_{\mc O})^{\rm T}\in \mc{C},$ there exists a unique solution $(w, y_{\mc O}, \dot w, \dot y_{\mc O})^{\rm T} $ to (\ref{maincont})-(\ref{initialcont}) with $(w, y_{\mc O}, \dot w, \dot y_{\mc O})^{\rm T}\in C([0,T]; \mc{C})$
and
\begin{eqnarray}
\nonumber \|(w, y_{\mc O}, \dot w, \dot y_{\mc O})^{\rm T}\|_{\mc{C}}\le C\left\{ \|(w^0, y_{\mc O}^0,  w^1,  y_{\mc O}^1)^{\rm T}\|_{\mc{C}}+ \|(M, {\bf{g}}_{\mc O})\|_{(L^2(\Omega))^{(m+2)}}\right\}.
\end{eqnarray}
\end{prop}

Our main exact controllability theorem is the following:
\vspace{0.1in}

\begin{thm}  \label{regularity} Let $T> \tau$ where
 \begin{eqnarray}\label{tau}\tau:=2L \left[\mathop {\min }\limits_{i=1,3,\ldots, 2m+1}\left(\sqrt{\frac{K}{\alpha}},~ \sqrt{\frac{\rho_i}{E_i}}\right)\right]^{-1}.
\end{eqnarray} For sufficiently small $\|\tilde {\bf{ G}}_E\|$ and for any $(w^0, {y}^0_{\mc O}, w^1,  {y}^1_{\mc O})^{\rm T}\in \mc{C}$ there exists $(M(t), {\bf{g}}_{\mc O}(t)) \in (L^2(0,T))^{(m+2)}$ such that
$(w(T), {y}_{\mc O}(T), \dot w(T), \dot {y}(T))^{\rm T}=0.$
\end{thm}
 \vspace{0.1in}

Now consider
\begin{equation}\left\{ \begin{array}{l}
m\ddot z -\alpha \ddot z'' +  K   z'''' -   N^{\rm T} {\bf{h}}_E \left({\bf{ G}}_E \phi_E + \tilde{{\bf{ G}}}_E \dot \phi_E \right)' = 0 ~~~ {\rm{on}}~~ \Omega\times \mathbb{R}^+ \\
   {\bf{h}}_{\mc O} {\bf{p}}_{\mc O} {\ddot v}_{\mc O} -{\bf{h}}_{\mc O} {\bf{E}}_{\mc O}  {v}_{\mc O}'' + {\bf{B}}^{\rm T}  \left({\bf{ G}}_E \phi_E + \tilde{{\bf{ G}}}_E \dot \phi_E \right) = 0~~ {\rm{on}} ~~ \Omega\times \mathbb{R}^+ \\
{\mbox{where}}~~  {\bf{B}} { v}_{\mc O}={\bf{h}}_E \phi_E-{\bf{h}}_E  N z'
 \end{array} \right.
\label{mainhomo}
\end{equation}
 with either hinged-Neumann (h-N), or clamped-Dirichlet (c-D), or mixed-mixed (m-m) boundary conditions respectively
 \begin{subequations}\label{bccc}
\begin{empheq}[left={\phantomword[l]{}{ } \hspace{-0.2in}\empheqlbrace}]{align}
\label{bdry33}&z(0,t)=z''(0,t)= z(L,t)= z''(L,t) = 0, {v}_{\mc O}'(0,t) = {v}_{\mc O}'(L,t)= 0 & \text{(h-N)}\\
\label{bdry55}&z(0,t)=z'(0,t)= z(L,t)= z'(L,t) = 0 , {v}_{\mc O}(0,t) = {v}_{\mc O}(L,t)=0 & \text{(c-D)}\\
\label{bdry77}& z(0,t)=z'(0,t)= z(L,t)= z''(L,t) = 0 ,{v}_{\mc O}(0,t) =  {v}_{\mc O}'(L,t)=  0. & \text{(m-m)}
\end{empheq}
\end{subequations}
 The initial conditions for (\ref{mainhomo}) are
\begin{eqnarray}
z(x,0)=z^0(x)&,&  ~~\dot z(x,0)=z^1(x), ~~ {v}_{\mc O}(x,0)= {v}^0_{\mc O}, ~~ {\dot v}_{\mc O}(x,0)= {v}^1_{\mc O}.\label{initialhomo}
\end{eqnarray}
For convenience, let $\mc{S}$ be a set, and $f,g$ be nonnegative functions on $\mc{S}.$ We will write $f\asymp g$ if there exists $C>0$ such that
$$\frac{1}{C} f(\lambda) \le g(\lambda)\le Cf(\lambda), ~~ \forall \lambda \in \mc{S}. $$

    The results in Theorem \ref{regularity} are based upon the following observability and hidden regularity results:
\begin{thm}\label{observability}
Let $T > \tau.$  Then for sufficiently small $\|\tilde {\bf{G}}_E\|$ solutions of the problem  (\ref{mainhomo})- (\ref{initialhomo}) satisfy the following observability and hidden regularity estimates:
\begin{subequations}\label{obs}
\begin{empheq}[left={\phantomword[l]{}{ }  \empheqlbrace}]{align}
\label{ohinged} & \int_0^T \left(|z'''(L,t)|^2 + |v_{\mc O}''(L,t)|^2\right)~dt \asymp  \|(z^0, v_{\mc O}^0, z^1, v_{\mc O}^1)^{\rm T}\|^2_{\mc{H}}& \text{(h-N)}~~~ \\
\label{oclamped}& \int_0^T \left(|z''(L,t)|^2 + |v_{\mc O}'(L,t)|^2\right)~dt \asymp  \|(z^0, v_{\mc O}^0, z^1, v_{\mc O}^1)^{\rm T}\|^2_{\mc{H}} & \text{(c-D)}~~~\\
\label{omixed}& \int_0^T \left( |z'(L,t)|^2 + |v_{\mc O}(L,t)|^2\right) ~dt \asymp  \|(z^0, v_{\mc O}^0, z^1, v_{\mc O}^1)^{\rm T}\|^2_{\mc{H}_{-1}} & \text{(m-m)}~~~~
\end{empheq}
\end{subequations}
where $\mc{H}$ and $\mc{H}_{-1}$ are later defined in (\ref{semigroupdom}) and (\ref{H-1}), respectively.
\end{thm}

Our results are improvements on earlier results \cite{Rajaram-Hansen3}, \cite{Rajaram-Hansen2} in several regards.   Here, we consider  the general multilayer system.
 The restriction on the size of ${{\bf{ G}}}$ has been eliminated, there are no conditions on the wave speeds, and  the optimal
 control time (determined by characteristics) is obtained.

  Our overall methodology is to first obtain appropriate boundary observability estimates for the uncoupled system of equations.
  This part uses mainly known estimates for the wave equation together with observability results obtained in \cite{O-Hansen1}.
  Second,  we prove, based on  carefully picked complex multipliers, a uniqueness result (Lemma \ref{overtrivial}) for the over-determined eigensystem  of the coupled system without damping $\tilde{ \mathbf  G}=0$ consisting of the homogeneous boundary conditions together with zero observation.
This allows us to deduce (using Theorem 6.2  in \cite{Komornik-P})     observability of the  coupled system without damping.   Finally, we are able include the possibility of small damping by a perturbation argument.

We consider  three different sets of boundary conditions.     While the overall structure of the proofs are the same in each case,
 the spaces  that arise  are different and lead to some very different technical issues.   For example,  in the case of (h-N) boundary conditions,   the system is well-posed with respect to  a higher-order energy defined by an extra derivative applied to each variable.   This allows us to obtain (similar to \cite{Komornik}, \cite{Lagnese-Lions}, \cite{LT1}) an observability result in a correspondingly smooth space, which  is equivalent to controllability in the natural energy space.
This approach fails in the case of (m-m) boundary conditions, where instead, we  obtain an observability result for  weaker solutions in which certain orthogonality conditions arise  (see Lemma \ref{lem*}).     In the case of (c-D) boundary conditions we obtain an observability result in the standard energy space, which in turn corresponds to an exact controllability result in a weaker space involving a quotient  $\mathrm M$ in the velocity component  of the transverse displacement in (\ref{maincont}).
The quotient $\mathrm M$  can not be eliminated if $L^2(0,T)$ controls are used.   This is due to
  orthogonality conditions   on   the range of the operator  $\mc L\phi = m \phi  - \alpha \phi''$  on the domain $H^2_0(\Omega)$   which must be imposed
  in the  transpositional solution.       (See Section \ref{Results-other} for details.)        In fact, a quotient space analogous to $\mathrm M$ was found in the velocity component of the optimal controls for  boundary control
  of the Kirchhoff plate  with clamped boundary conditions, \cite{LT2}.         Related optimal controllability and observability results   for the Rayleigh beam  are described in   \cite{O-Hansen1}.

   All of the  controllability results in this paper are optimal in the sense that the space of exact controllability matches the optimal regularity space for $L^2(0,T)$ boundary controls.   Moreover,   as mentioned above, the quotient $\mathrm M$  in (\ref{control2})  can not be eliminated from the control space if $L^2(0,T)$ controls are used.      On the other hand,  the quotients that occur in the second and fourth components of the control space (\ref{control1}) are perhaps inessential in  that  they   arise as a consequence of orthogonality constraints imposed for convenience in the homogeneous solutions  (see   (\ref{ohinged2})) which are used in the definition of transpositional solution (see Definition 5.1).   In this case solutions in    (\ref{control1})  are defined up to uniform translational motion in each layer.

   This paper is organized as follows.
In Section \ref{C2S2-SemFor} we prove regularity results for the homogeneous system using  semigroup theory.   In Section \ref{C2S2A-Extens} we characterize the weaker observability  space for the case of  (m-m) boundary conditions.     In Section \ref{C2S3-Observ} we prove  the key uniqueness result Lemma \ref{xyz}  and main observability result Theorem \ref{observability}.       In Section \ref{C2S4-Control} we define transpositional  solutions of the control problem and prove our main controllability result Theorem \ref{regularity}.

\section{Semigroup formulation}\label{C2S2-SemFor} Let  $$U=:(u,{\bf{u}})^{\rm T}=(z, {v}_{\mc O})^{\rm T}, ~~~~V:= (v, {\bf{v}})^{\rm T}=(\dot z,  {\dot v}_{\mc O})^{\rm T}, ~~{\rm{and}} ~~ Y:=(U,V)^{\rm T}.$$ Let $\mc L\varphi=m \varphi - \alpha \varphi''$. From the Lax-Milgram theorem $\mc L: H^1_0(\Omega)\to H^{-1}(\Omega)$ is an isomorphism which remains isomorphic from $H^2(\Omega)\cap H^1_0(\Omega)$ to $L^2(\Omega).$\\
Then {(\ref{mainhomo})-(\ref{initialhomo})} can be written as
\begin{eqnarray} &&\frac{dY}{dt}=\mc{A} Y:=\left( {\begin{array}{*{20}c}
   {0} & {I }  \\
   {-A_1 } & {A_2}    \\
\end{array}} \right) \left( \begin{array}{l}
 U \\
 V \\
 \end{array} \right),  Y(0)=(U(0),V(0))^{\rm T}=(z^0,  {v}^0_{\mc O}, z^1,  {v}^1_{\mc O} )^{\rm T}\quad\quad
\label{semigroupfor1}\end{eqnarray}
 where
\begin{eqnarray}A_1 U:=\left( {\begin{array}{*{20}c}
   \mc L^{-1}\left( Ku''''- N^{\rm T} {\bf{h}}_E {\bf{G}}_E ({\bf{h}}_E^{-1} {\bf{B}} {\bf{u}}'+  N u'') \right)  \\
   { {\bf{h}}^{-1}_{\mc O} {\bf{p}}^{-1}_{\mc O} } \left(-{ {\bf{h}}_{\mc O} {\bf{E}}_{\mc O} }\textbf{u}''+{\bf{B}}^{\rm T} {\bf{G}}_E( {\bf{h}}^{-1}_E {\bf{B}} {\bf{u}}+ N u')\right) \\
\end{array}} \right),\label{A1}\end{eqnarray}
\begin{eqnarray}\nonumber A_2 V:=\left( {\begin{array}{*{20}c}
  \mc  L^{-1}\left( N^{\rm T} {\bf{h}}_E \tilde {\bf{G}}_E ({\bf{h}}_E^{-1} {\bf{B}} {\bf{v}}'+  N v'') \right)  \\
   { {\bf{h}}^{-1}_{\mc O} {\bf{p}}^{-1}_{\mc O} } \left(-{\bf{B}}^{\rm T} \tilde{\bf{G}}_E( {\bf{h}}^{-1}_E {\bf{B}} {\bf{v}}+ N v')\right) \\
\end{array}} \right).
\end{eqnarray}

Let $\left<u,v\right>_{\Omega}=\int_{\Omega} u \cdot \overline{v}~ dx$ where $u$ and $v$ may be scalar or vector valued. Define the bilinear forms $a$ and $c$ by
\begin{eqnarray}
\nonumber c(z, v_{\mc O}; \hat{z}, \hat{v}_{\mc O} )&=& m\left<z,\hat{z}\right>_{\Omega}+\alpha\left<z',\hat{z}'\right>_{\Omega}+ \left< {\bf{h}}_{\mc O} {\bf{p}}_{\mc O} { v}_{\mc O},   {\hat{ v}}_{\mc O}\right>_{\Omega},\\
\nonumber  a(z, v_{\mc O}; \hat{z}, \hat{v}_{\mc O} )&=& K  \left<z'',\hat{z}''\right>_{\Omega}+ \left< {\bf{h}}_{\mc O} {\bf{E}}_{\mc O} {v}_{\mc O}',   { \hat{v}}_{\mc O}'\right>_{\Omega}+   \left<{\bf{G}}_E {\bf{h}}_E \phi_E,{\hat{\phi}}_E\right>_{\Omega}\\
\nonumber & =&K  \left<z'',\hat{z}''\right>_{\Omega}+ \left< {\bf{h}}_{\mc O} {\bf{E}}_{\mc O} {v}_{\mc O}',   { \hat{v}}_{\mc O}'\right>_{\Omega}\\
\label{forms} && \quad+  \left<{\bf{G}}_E {\bf{h}}_E^{-1}\left( {\bf{B}} { v}_{\mc O}+N z'\right),\left( {\bf{B}} {\hat{v}}_{\mc O}+N \hat z'\right)\right>_{\Omega}.
\end{eqnarray}
The ``higher order'' and natural  energies of the beam are respectively given by
 \begin{subequations}\label{energy}
\begin{empheq}[left={\phantomword[r]{}{\mc{E}(t) = }  \empheqlbrace}]{align}
\label{energyhN}& \frac{1}{2}\left( a(z', v_{\mc O}')+c(\dot z', \dot v_{\mc O}')\right) & \text{(h-N)} \\
\label{energycD} & \frac{1}{2}\left(a(z,v_0) + c(\dot z, \dot v_{\mc O})\right) & \text{(c,D), (m-m)},
\end{empheq}
\end{subequations}
where $a(\cdot), c(\cdot)$ are the quadratic forms that agree with $a(\cdot, \cdot), c(\cdot, \cdot)$ on the diagonal. Define the energy inner products corresponding to each set of boundary conditions by
 \begin{subequations}\label{inner}
\begin{empheq}[left={\phantomword[r]{}{\left<Y,\widehat{Y}\right>_{\mc{H}}=}  \empheqlbrace}]{align}
\label{innerhN}& a(U';\widehat{U}')+c(V';\widehat{V}'). & \text{(h-N)}\\
\label{innercD} & a(U;\widehat{U})+c(V;\widehat{V}) & \text{(c-D)}, \text{(m-m)}.
\end{empheq}
\end{subequations}
Corresponding to each case, define the Hilbert spaces
 \begin{subequations}\label{semigroupdom}
\begin{empheq}[left={\phantomword[r]{}{\mc{H} = }  \empheqlbrace}]{align}
\label{ohinged2}& H^3_*(\Omega)\times \left(H^2_\perp(\Omega)\right)^{(m+1)}\times \left(H^2(\Omega)\cap H^1_0(\Omega)\right)\times (H^1_\perp(\Omega))^{(m+1)} & \text{(h-N)}\\
\label{oclamped2} &  H^2_0(\Omega)\times \left(H^1_0(\Omega)\right)^{(m+1)}\times H^1_0(\Omega)\times (L^2(\Omega))^{(m+1)} & \text{(c-D)} \\
\label{omixed2}& H^2_\#(\Omega)\times \left(H^1_\dag(\Omega)\right)^{(m+1)}\times \left(H^1_0(\Omega)\right)\times (L^2(\Omega))^{(m+1)} & \text{(m-m)}
\end{empheq}
\end{subequations}
where $H^2_\#(\Omega)$ and $H^1_\dag(\Omega)$ are defined in (\ref{spaces}) and
\begin{eqnarray}
\nonumber && H^3_*(\Omega):=\{ u \in H^3(\Omega) \cap H_{0}^1(\Omega)~:~ u''(0)=u''(L)=0\}\\
\nonumber && H^1_\perp(\Omega):=\{u\in H^1(\Omega)~:~\int_{\Omega} u~dx=0\}\\
\nonumber && H^2_\perp(\Omega):=\{u\in H^2(\Omega)\cap H^1_\perp(\Omega)~:~~ u'(0)=u'(L)=0\}.
\end{eqnarray}
Define $ \mc{D}(\mc{A})$ by
 \begin{subequations}\label{semigroupdom2}
\begin{empheq}[left={\phantomword[r]{}{\mc{D}( \mc{A})=  }  \empheqlbrace}~~]{align}
\nonumber & \left(H^4(\Omega)\cap H^3_*(\Omega)\right)\times \left(H^3(\Omega)\cap H^2_\perp(\Omega)\right)^{(m+1)}\times H^3_*(\Omega)\times (H^2_\perp(\Omega))^{(m+1)} & \text{(h-N)}\\
\nonumber & \left(H^3(\Omega)\cap H^2_0(\Omega)\right)\times \left(H^2(\Omega)\cap H^1_0(\Omega)\right)^{(m+1)}\times H^2_0(\Omega)\times (H^1_0(\Omega))^{(m+1)}  & \text{(c-D)} \\
\nonumber & H^3_\#(\Omega) \times \left(H^2_\dag(\Omega)\right)^{(m+1)}\times H^2_\#(\Omega)\times (H^1_\dag(\Omega))^{(m+1)} & \text{(m-m)}
\end{empheq}
\end{subequations}
where
\begin{eqnarray}\nonumber && H^3_\#(\Omega):=\{ u \in H^2_\#(\Omega) ~:~ u''(L)=0\},\\
\nonumber && H^2_\dag(\Omega):=\{u\in H^2(\Omega)\cap H^1_\dag(\Omega)~:~u'(L)=0\}.
\end{eqnarray}

\vspace{0.1in}
\begin{Lem} {\label{densdefined}}The operator $\mc{A}: \mc{D}(\mc{A})\subset \mc{H} \to \mc{H}$ is densely defined.
\end{Lem}
\vspace{0.1in}

\textbf{Proof:} The density is obvious. However, in the case of hinged-Neumann boundary conditions (h-N), it is not obvious that the orthogonality constraint in the definition of $\mc H$ is invariant with respect to $\mc A,$ i.e., that $Y\in \mc D(\mc A)$ implies $\mc A Y\in \mc H.$ To verify this, let $Y=(u, {\bf{u}},v,{\bf{v}})^{\rm T} \in \mc D (\mc A).$ Then
 $$(u, {\bf{u}},v,{\bf{v}})^{\rm T}\in  \left(H^4(\Omega)\cap H^3_*(\Omega)\right)\times \left(H^3(\Omega)\cap H^2_\perp(\Omega)\right)^{(m+1)}\times H^3_*(\Omega)\times (H^2_\perp(\Omega))^{(m+1)}.$$ From (\ref{semigroupfor1}), $\mc{A}Y= \left( \begin{array}{c}
 V\\
 0 \\
 \end{array} \right) +  \left( \begin{array}{c}
 0 \\
 -A_1 U + A_2 V \\
 \end{array} \right).$ Since $v\in H^3_*(\Omega) $ and ${\bf v}\in (H^2_\perp(\Omega))^{(m+1)},$ $\left( \begin{array}{c}
 V\\
 0 \\
 \end{array} \right)\in \mc{H}. $ Explicitly,  $-A_1 U + A_2 V $ is
\begin{eqnarray}\label{nut}  \left( \begin{array}{c}
\mc L^{-1}\left( -Ku''''+ N^{\rm T} {\bf{h}}_E \left[{\bf{G}}_E ({\bf{h}}_E^{-1} {\bf{B}} {\bf{u}}'+  N u'') +   \tilde {\bf{G}}_E ({\bf{h}}_E^{-1} {\bf{B}} {\bf{v}}'+  N v'') \right] \right)   \\
 { {\bf{h}}^{-1}_{\mc O} {\bf{p}}^{-1}_{\mc O} } \left({ {\bf{h}}_{\mc O} {\bf{E}}_{\mc O} }\textbf{u}''-{\bf{B}}^{\rm T}\left[ {\bf{G}}_E( {\bf{h}}^{-1}_E {\bf{B}} {\bf{u}}+ N u') - \tilde{\bf{G}}_E( {\bf{h}}^{-1}_E {\bf{B}} {\bf{v}}+ N v')\right]\right)
 \end{array} \right).~~~~
 \end{eqnarray}
 The first entry of (\ref{nut}) is in $\left(H^2(\Omega)\cap H^1_0(\Omega)\right)$ since $\mc L^{-1}$ maps $L^2(\Omega)$ to \\$\left(H^2(\Omega)\cap H^1_0(\Omega)\right).$ Lastly, the second entry of (\ref{nut}) is in $(H^1_\perp(\Omega))^{(m+1)} $ since the application of the (h-N) boundary conditions implies $\int_\Omega u' ~dx = \int_\Omega{\bf u}'' ~dx =0.$ Furthermore, since $Y\in \mc{D}(\mc{A}),$ $\int_\Omega {\bf u} ~dx = \int_\Omega{\bf v} ~dx =0,$ it follows that
 \begin{eqnarray}
 \nonumber  \int_\Omega { {\bf{h}}^{-1}_{\mc O} {\bf{p}}^{-1}_{\mc O} }{\bf{B}}^{\rm T} {\bf{G}}_E {\bf{h}}^{-1}_E {\bf{B}} {\bf{u}}~ dx = \int_\Omega { {\bf{h}}^{-1}_{\mc O} {\bf{p}}^{-1}_{\mc O} }{\bf{B}}^{\rm T} {\tilde{\bf{G}}}_E {\bf{h}}^{-1}_E {\bf{B}} {\bf{v}}~ dx= 0.~~\square
 \end{eqnarray}
 \vspace{0.1in}

 \begin{Lem} \label{skew-adjoint} The infinitesimal generator $\mc{A}$ for each set of boundary conditions is dissipative, and moreover it satisfies
  \begin{subequations}\label{dissipation}
\begin{empheq}[left={\phantomword[c]{}{ {\mbox{Re}}\left<\mc{A}Y, Y\right>_{\mc{H}}=}\quad\quad\quad \empheqlbrace}~]{align}
&- \left<  \tilde {\bf{G}}_E \Theta', {\bf{h}}_E^{-1}\Theta' \right>_{\Omega}\le 0,& ~  \text{(h-N)}\\
&- \left<  \tilde {\bf{G}}_E \Theta, {\bf{h}}_E^{-1}\Theta\right>_{\Omega}\le 0,&~  \text{(c-D)}, \text{(m-m)}
\end{empheq}
\end{subequations}
for all $Y=(u, {\bf{u}},v,{\bf{v}})^{\rm T}\in \mc D(\mc A)$ where $\Theta=\left( {\bf{B}} {\bf{v}} + {\bf{h}}_E    N v' \right).$
 \end{Lem}
\vspace{0.1in}

\textbf{Proof:}
It is easy to show that
 $\mc{A}$ is dissipative on $\mc{H}$ for each set of boundary conditions. For example, consider the (h-N) boundary conditions:
  \begin{eqnarray}
\nonumber \left<\mc{A}Y, Y\right>_{\mc{H}} &=&   \left\{-K \left<u''',  v'''\right>_{\Omega} +K\left<{ v}''', u'''\right>_{\Omega}\right\}+\left\{-\left<{ {\bf{h}}_{\mc O} {\bf{E}}_{\mc O} }\textbf{u}'', {\bf{ v}}''\right>_{\Omega}+ \left<{\bf{h}}_{\mc O} {\bf{E}}_{\mc O} {\bf{ v}}'',  {\bf{u}''}\right>_{\Omega}\right\}\\
\nonumber   && +\left\{- \left< {\bf{G}}_E\left( {\bf{B}} {\bf{u}}' + {\bf{h}}_E   N u''\right), {\bf{h}}^{-1}_E \left({\bf{B}} {\bf{ v}'} + {\bf{h}}_E^{-1}N v''\right)\right>_{\Omega} \right. \\
 \nonumber &&~~~\left.+ \left<{\bf{G}}_E  ({\bf{B}} {\bf{ v}}'+  {\bf{h}}_E N v''), {\bf{h}}_E^{-1} ({\bf{B}} {\bf{u}}'+ {\bf{h}}_E N u'')\right>_{\Omega}\right\}\\
\nonumber &&- \left<  \tilde {\bf{G}}_E \left( {\bf{B}} {\bf{ v}}' + {\bf{h}}_E    N \bar v'' \right),  {\bf{h}}_E^{-1}\left({\bf{B}} {\bf{ v}'} + {\bf{h}}_E N  v''\right)\right>_{\Omega} \\
\nonumber &=& -2i~ {\rm Im} \left( K\left<u''', v'''\right>_{\Omega}\right) - 2i~ {\rm Im}\left(\left< {\bf{h}}_{\mc O} {\bf{E}}_{\mc O}{\bf{u}''}, {\bf{ v}}''\right>_{\Omega}\right)\\
 \nonumber && - 2i~ {\rm Im} \left< {\bf{G}}_E\left( {\bf{B}} {\bf{u}}' + {\bf{h}}_E   N u''\right), {\bf{h}}^{-1}_E \left({\bf{B}} {\bf{ v}'} + {\bf{h}}_E^{-1}N v''\right)\right>_{\Omega}  - \left<  \tilde {\bf{G}}_E \Theta', {\bf{h}}_E^{-1}\Theta' \right>_{\Omega}.
\end{eqnarray}
Therefore (\ref{dissipation}) follows. $\square$
     \vspace{0.1in}

\begin{Lem} \label{surjective} $I-\mc{A}: \mc{D}(\mc{A})\to \mc{H}$ is surjective.
\end{Lem}
\vspace{0.1in}

\textbf{Proof:} We prove the lemma for only (h-N) boundary conditions since the proofs for other boundary conditions are similar. Let $C$ denote a generic constant in the following calculations, and define $|u|_s=\|u\|_{H^s(\Omega)}, ~|{\bf{u}}|_s=\|{\bf{u}}\|_{(H^s(\Omega))^{(m+1)}}.$ Let $Y_1=(u_1,{\bf{u}_1}, v_1, {\bf{v}_1})^{\rm T}.$  For given $Y_2=(u_2,{\bf{u}_2}, v_2, {\bf{v}_2})^{\rm T} \in  \mc{H} $ we want to prove the solvability of the system $(I-\mc{A})Y_1=Y_2$  in $\mc{D}({\mc{A}}):$ 
\begin{eqnarray}
\nonumber  Ku_1''''- N^{\rm T} {\bf{h}}_E \left({\bf{G}}_E ({\bf{h}}_E^{-1} {\bf{B}} {\bf{u}_1}'+  N u_1'') +  \tilde {\bf{G}}_E ({\bf{h}}_E^{-1} {\bf{B}} {\bf{v}_1}'+  N v_1'')\right)&&= \mc Lv_2-\mc Lv_1 \\
\nonumber -{ {\bf{h}}_{\mc O} {\bf{E}}_{\mc O} }\textbf{u}_1''+{\bf{B}}^{\rm T} \left({\bf{G}}_E( {\bf{h}}^{-1}_E {\bf{B}} {\bf{u}_1}+ N u_1')+ \tilde{\bf{G}}_E( {\bf{h}}^{-1}_E {\bf{B}} {\bf{v}_1}+ N v_1') \right)&&=  {\bf{p}}_{\mc O} {\bf{h}}_{\mc O} \left({\bf{v}_2}-{\bf{v}_1}\right)\\
 \nonumber  u_1-v_1 && =u_2 \\
\label{rangecon}  {\bf{u}_1} - {\bf{v}_1}&&  = {\bf{u}_2}.
\end{eqnarray}
Differentiating the second equation in (\ref{rangecon}) yields
   \begin{eqnarray}
\nonumber  Ku_1''''- N^{\rm T} {\bf{h}}_E \left({\bf{G}}_E ({\bf{h}}_E^{-1} {\bf{B}} {\bf{u}_1}'+  N u_1'') +  \tilde {\bf{G}}_E ({\bf{h}}_E^{-1} {\bf{B}} {\bf{v}_1}'+  N v_1'')\right)&&= \mc Lv_2-\mc Lv_1 \\
\nonumber -{ {\bf{h}}_{\mc O} {\bf{E}}_{\mc O} }\textbf{u}_1'''+{\bf{B}}^{\rm T} \left({\bf{G}}_E( {\bf{h}}^{-1}_E {\bf{B}} {\bf{u}_1'}+ N u_1'')+ \tilde{\bf{G}}_E( {\bf{h}}^{-1}_E {\bf{B}} {\bf{v}_1'}+ N v_1'') \right)&&=  {\bf{p}}_{\mc O} {\bf{h}}_{\mc O} \left({\bf{v}_2'}-{\bf{v}_1'}\right)\\
 \nonumber  u_1-v_1 && =u_2 \\
\label{rangecon1}  {\bf{u}_1} - {\bf{v}_1}&&  = {\bf{u}_2}.
\end{eqnarray}
We eliminate the functions $v_1, {\bf{v}_1}$ from the last two equations in (\ref{rangecon1}). Then, we multiply the first equation $u_1''''$ and the second  by ${\bf{u}_1}''',$ and integrate by parts on $\Omega,$ using boundary conditions for $\mc{D}(\mc{A}),$ and then we eventually use Holder's inequality to obtain the following estimate:
\begin{eqnarray}
\nonumber |u_1|_4 \le &&~ C\left( |u_1|_2+|{\bf{u}_1}|_1+|u_2|_2+|v_2|_2 + |{\bf{u}_2}|_1\right)\\
\nonumber |\textbf{u}_1|_3 \le && ~C\left(|u_1|_2+|{\bf{u}_1}|_1+|u_2|_2 + |{\bf{u}_2}|_2+|{\bf{v}_2}|_1\right)\\
\nonumber |v_1|_3 \le &&~ C\left( |u_1|_3+|u_2|_3 \right)\\
\label{rangecon2} |{\bf{v}}_1|_2 \le &&~ C\left( |{\bf{u}}_1|_2+|{\bf{u}}_2|_2\right).
\end{eqnarray}
The next step is to absorb the lower order terms in (\ref{rangecon2}) to get
\begin{equation}
\label{rangecon3} |u_1|_4  + |\textbf{u}_1|_3  + |v_1|_3 + |{\bf{v}}_1|_2\le C\left(|u_2|_3  + |\textbf{u}_2|_2  + |v_2|_2 + |{\bf{v}}_2|_1\right).
\end{equation}
  We apply a standard compactness-uniqueness argument: now suppose contrarily that the inequality (\ref{rangecon3}) does not hold. Then there exists a sequence $Y_{2n}:=\{(u_{2n}, \textbf{u}_{2n}, v_{2n} , {\bf{v}}_{2n})^{\rm T}\}_{n=1}^{\infty}$ such that
\begin{equation}\|Y_{2n}\|_{\mc{H}} \mathop  \to 0, \quad{\text{and}}\quad |u_{1n}|_4  + |\textbf{u}_{1n}|_3  + |v_{1n}|_3 + |{\bf{v}}_{1n}|_2=1.\label{rangecon4}\end{equation} and
From (\ref{rangecon4}) we can extract a subsequence, still denoted \\$Y_{1n}:=\{[u_{1n}, \textbf{u}_{1n}, v_{1n} , {\bf{v}}_{1n}]^{\rm T}\}_{n=1}^{\infty}$  such that $Y_{1n}$ converges to $Y_1:=(u_1, {{\bf u}_1}, v_1, {{\bf v}_1})$ weakly in $ H^4(\Omega)\times \left(H^3(\Omega)\right)^{(m+1)}\times H^3(\Omega)\times (H^2(\Omega))^{(m+1)}:=\mc{W}.$
If we consider the solution of (\ref{rangecon}) with  $Y_{1n}=Y_{1n}(Y_{2n}),$  then it follows from (\ref{rangecon2}) that
\begin{eqnarray}
\nonumber  |u_{1n}-u_{1m}|_4 &\le&  C \left(|u_{1n}-u_{1m}|_2+|{\bf{u}_{1n}}-{\bf{u}_{1m}}|_1+|u_{2n}-u_{2m}|_2\right.\\
\nonumber && \left.\quad\quad+|v_{2n}-v_{2m}|_2 + |{\bf{u}_{2n}}-{\bf{u}_{2m}}|_1\right) \\
\nonumber  |\textbf{u}_{1n}-\textbf{u}_{1m}|_3 &\le& C\left(|u_{1n}-u_{1m}|_2+|{\bf{u}_{1n}}-{\bf{u}_{1m}}|_1+|u_{2n}-u_{2m}|_2\right.\\
 \nonumber &&\quad\quad\left.+ |{\bf{u}_{2n}}-{\bf{u}_{2m}}|_2+|{\bf{v}_{2n}}-{\bf{v}_{2m}}|_1\right)\\
\nonumber |v_{1n}-v_{1m}|_3 &\le&  C\left( |u_{1n}-u_{1m}|_3+|u_{2n}-u_{2m}|_3\right) \\
\label{rangecon8}  |{\bf{v}}_{1n}-{\bf{v}}_{1m}|_2 &\le& C\left( |{\bf{u}}_{1n}-{\bf{u}}_{1m}|_2+|{\bf{u}}_{2n}-{\bf{u}}_{2m}|_2\right).
\nonumber
\end{eqnarray}
Thus, by the Sobolev's compact embedding theorem we get
$$|u_{1n}-u_{1m}|_4, ~|\textbf{u}_{1n}-\textbf{u}_{1m}|_3,~ |v_{1n}-v_{1m}|_3,~ |{\bf{v}}_{1n}-{\bf{v}}_{1m}|_2 \to 0, $$ as $n,m\to\infty.$ This implies that  $Y_{1n}$ actually converges to $Y_1$ strongly in $\mc{W}.$
On the other hand, the system (\ref{rangecon}) with $Y_2=(0,{\bf{0}}, 0, {\bf{0}})^{\rm T}$, see (\ref{rangecon4}), has only a trivial solution since the system (\ref{semigroupfor1}) is dissipative by (\ref{dissipation}). This contradicts with (\ref{rangecon4}) and therefore (\ref{rangecon3}) holds. Hence $Y_1\in \mc{D}(\mc{A})$ and the claim of the theorem is proved.
\vspace{0.1in}

\begin{thm} \label{eigens}$\mc{A}: \mc{D}( \mc{A}) \to \mc{H}$ is the infinitesimal generator of a $C_0-$semigroup of contractions.  Moreover, the spectrum of $\mc{A}$  only consists of isolated non-zero eigenvalues $\{\gamma_n\}_{n=1}^\infty, $ and $|\gamma_n^{\pm}|\to \infty $ as $n\to\infty.$
\end{thm}
\vspace{0.1in}

\textbf{Proof:} The proof  of the first part follows from the  L\"{u}mer-Phillips theorem \cite{Pazy}  using Lemma \ref{densdefined}, \ref{skew-adjoint} and \ref{surjective}. Since $(\mc I-\mc A)^{-1}$ is compact, the spectrum of $\mc A$ only consists of eigenvalues. A simple proof that  $0\in\rho(\mc A)$  for the (h-N) case ($m=1$) is given in \cite{Rajaram-Hansen3}. The same proof applies for any positive integer $m$ and also the boundary conditions (c-D) and (m-m).  Hence the claim of the theorem follows. $\square$
\vspace{0.1in}

 \begin{Cor} \label{SAdjoint} The operator $\mc {A}^* : \mc D(\mc A) =\mc D(\mc A^*)\to \mc H $ is the generator of a $C_0-$contraction semigroup. Moreover,
$$\left[\mc{A}(\tilde{{\bf{G}}}_E)\right]^*=-\mc{A}(-\tilde{{\bf{G}}}_E)), ~~\mbox{on} ~~ \mc{D}(\mc{A})= \mc{D}(\mc{A}^*)$$
where $\mc{A}(\tilde{{\bf{G}}}_E))$ denotes the dependence of $\mc{A}$ on the parameter $\tilde{{\bf{G}}}_E.$
\end{Cor}
\vspace{0.1in}

\textbf{Proof:}  A straightforward (but lengthy) calculation shows that  $\left[\mc{A}(\tilde{{\bf{G}}}_E)\right]^*=-\mc{A}(-\tilde{{\bf{G}}}_E)$ on $\mc{D}(\mc{A})$ for each of the sets of boundary conditions considered. Moreover $-\mc{A}(-\tilde{{\bf{G}}}_E)$ is dissipative by (\ref{dissipation}). Thus the proof of Lemma \ref{surjective} remains valid with $-\mc{A}(-\tilde{{\bf{G}}}_E)$ in place of $\mc A.$ Since $\mc I + \mc{A}(-\tilde{{\bf{G}}}_E): \mc D(\mc A) \to \mc H$ is bijective, $\mc D(\mc A^*)$ can be no larger than $\mc D(\mc A).$ Thus, $\mc D(\mc A^*)=\mc D(\mc A). $ It follows from the corollary of L\"{u}mer-Phillips theorem (\cite{Pazy}, Chap I) that $\mc A^*$ generates a contraction semigroup. $\square$
\vspace{0.1in}

Let $\mc H_{-1}$ be the dual space of $\mc H_1:= \mc D(\mc A)$ pivoted with respect to $\mc H.$ Then we have the following dense and compact embeddings
$$\mc H_1\subset \mc H\subset \mc H_{-1}.$$
By Proposition 2.10.3 in \cite{Weiss-Tucsnak}, the operator $\mc{A}: \mc H_1 \to \mc{H}$ has a unique extension   $\tilde{\mc A}: \mc H \to \mc{H}_{-1}$ defined by
 \begin{eqnarray} \left<\tilde {\mc{A}} Y, Z\right>:= \left< Y,  \mc{A}^* Z\right>_{\mc{H}}, ~~  \forall  ~Z \in \mc H_1, Y\in \mc{H}.
\label{definofH-1}
\end{eqnarray}
By Proposition 2.10.4 in \cite{Weiss-Tucsnak},   $\tilde{\mc{A}}$ is the generator of a  $C_0-$semigroup  $\{ e^{\tilde{ \mc{A}}t} \}_{t\ge 0}$ on $\mc H_{-1}$,which is similar to  $\{e^{\mc{A}t}\}_{t\ge 0}$.      Thus we have the following.
\vspace{0.1in}

 \begin{Cor}\label{ext} The semigroup $\{e^{\mc A t}\}_{t\ge 0}$ with the generator $\mc{A}: \mc H_1 \to \mc{H}$ has a unique extension to a  contraction semigroup    $\{e^{\tilde{\mc A }t}\}_{t\ge 0}$    on $\mc{H}_{-1}$  with the generator $\tilde{\mc{A}}: \mc H \to  \mc{H}_{-1}. $
\end{Cor}
\vspace{0.1in}

\section{ Characterization of the space $\mc H_{-1}$  in undamped case}
\label{C2S2A-Extens}
\vspace{0.1in}
In particular,  we are interested in a characterization of the space $\mc H_{-1}$  for the (m-m) boundary conditions. Define spaces $\mc{X}_2$, $\mc{X}_1$, $\mc{X}$    by
\begin{subequations}
\begin{empheq}[left={\phantomword[r]{}{ \mc{X}_2} = \empheqlbrace}]{align}
\nonumber & \left(H^4(\Omega)\cap H^3_*(\Omega)\right) \times ( H^3(\Omega)\cap H^2_\perp(\Omega))^{(m+1)}   & \text{(h-N)}\\
\nonumber &  \left(H^3(\Omega)\cap H^2_0(\Omega)\right) \times  (H^2(\Omega)\cap H^1_0(\Omega))^{(m+1)} & \text{(c-D)}\\
\nonumber & H^3_\#(\Omega) \times ( H^2_\dag(\Omega))^{(m+1)} & \text{(m-m)}
\end{empheq}
\begin{empheq}[left={\phantomword[r]{}{ \mc{X}_1} = \empheqlbrace}]{align}
\nonumber & H^3_*(\Omega) \times (H^2_\perp(\Omega))^{(m+1)}   & \text{(h-N)}\\
\nonumber & H^2_0(\Omega)\times (H^1_0(\Omega))^{(m+1)} & \text{(c-D)}\\
\nonumber & H^2_\#(\Omega) \times ( H^1_\dag(\Omega))^{(m+1)} & \text{(m-m)}
\end{empheq}
\begin{empheq}[left={\phantomword[r]{}{ \mc{X}} = \empheqlbrace}]{align}
\nonumber & \left(H^2(\Omega) \cap H^1_0(\Omega)\right) \times (H^1_\perp(\Omega))^{(m+1)}   & \text{(h-N)}\\
\nonumber & H^1_0(\Omega)\times (L^2(\Omega))^{(m+1)} & \text{(c-D)}, \text{(m-m).}
\end{empheq}
\end{subequations}

Also define the  inner products
 \begin{subequations}
\begin{empheq}[left={\phantomword[r]{}{ \left<U,V\right>_{\mc{X}_1}} = \empheqlbrace}]{align}
\nonumber & a(U';V')   & \text{(h-N)}\\
\nonumber & a(U;V) & \text{(c-D), (m-m)},
\end{empheq}
\end{subequations}
where   $U=(u,{\bf{u}})^{\rm T}, V=(v,{\bf{v}})^{\rm T} $  and the bilinear  form $a$ is defined in (\ref{forms});
  \begin{eqnarray}
  \label{eq2003}
 &  \left<U,V\right>_{\mc{X}}=\left\{
 \begin{array}{l}
  c(U';V')=m\left<u',v'\right>_{\Omega}+\alpha\left<u'',v''\right>_{\Omega}+ \left< {\bf{h}}_{\mc O} {\bf{p}}_{\mc O} {\bf{u}}',  {\bf{v}}' \right>_{\Omega}\\
\quad\quad\quad\quad=-\left<\mc L u,v''\right>_{\Omega}+ \left< {\bf{h}}_{\mc O} {\bf{p}}_{\mc O} {\bf{u}}',   {\bf{v}}'\right>_{\Omega}, \quad\quad\text{(h-N)}\\
 c(U;V)=m\left<u,v\right>_{\Omega}+\alpha\left<u',v'\right>_{\Omega}+ \left< {\bf{h}}_{\mc O} {\bf{p}}_{\mc O} {\bf{u}},  {\bf{v}} \right>_{\Omega}\\
~\quad\quad\quad=-\left<\mc L u,v\right>_{\Omega}+ \left< {\bf{h}}_{\mc O} {\bf{p}}_{\mc O} {\bf{u}},   {\bf{v}}\right>_{\Omega}, \quad\quad\text{(c-D), \text{(m-m).}}
   \end{array} \right. &
\end{eqnarray}
Then
\begin{equation} \nonumber \mc{D}(\mc{A})  =  \mc{X}_2\times \mc{X}_1,\qquad \mc{H} =   \mc{X}_1\times \mc{X}
\end{equation}
and     the inner product for $\mc H$ can be written
\begin{equation}\nonumber \left<Y,\hat Y\right>_{\mc H} = \left<(U,V)^T, (\hat U, \hat V)^T\right>_{\mc H}= \left<U, \hat U\right>_{{\mc X}_1} +  \left<V, \hat V\right>_{{\mc X}}.
\end{equation}

Let $A_1$ be the operator on $\mc{X}_1$ defined by (\ref{A1}).
For each of the sets of  boundary conditions (h-N), (m-m) or (c-D), a simple calculation establishes the following identity:
\begin{equation}
 \left<A_1 U,V\right>_{\mc{X}}  =  \left<U,V\right>_{\mc{X}_1}\qquad \forall \, U,V\in \mc{X}_2. \label{eq2001}
 \end{equation}
 For instance,  in the (h-N) case,
\begin{eqnarray*}\nonumber \left<A_1 U,V\right>_{\mc{X}}&&=\left<\left( {\begin{array}{*{20}c}
   \mc L^{-1}\left( Ku''''- N^{\rm T} {{\bf{h}}_E {\bf{G}}_E} \phi_E' \right)  \\
   { {\bf{h}}^{-1}_{\mc O} {\bf{p}}^{-1}_{\mc O} } \left(-{ {\bf{h}}_{\mc O} {\bf{E}}_{\mc O} }\textbf{u}''+{\bf{B}}^{\rm T} {\bf{G}}_E \phi_E\right) \\
\end{array}} \right),V \right>_{\mc{X}}\\
\nonumber && =\left<-Ku''''+ N^{\rm T} {\bf{h}}_E {\bf{G}}_E \phi_E', v''\right>_{\Omega}+ \left<-{ {\bf{h}}_{\mc O} {\bf{E}}_{\mc O} }\textbf{u}'''+{\bf{B}}^{\rm T} {\bf{G}}_E \phi_E', {\bf{v}}'\right>_{\Omega}\\
\nonumber && = K\left<u''', v'''\right>_{\Omega}+ \left<{ {\bf{h}}_{\mc O} {\bf{E}}_{\mc O} }\textbf{u}'', {\bf{v}}''\right>_{\Omega}+  \left<  {\bf{G}}_E \phi_E', {\bf{h}}_E N v'' + {\bf{B}} {\bf{v}'}\right>_{\Omega}  \\
 && = K\left<u''', v'''\right>_{\Omega}+ \left<{ {\bf{h}}_{\mc O} {\bf{E}}_{\mc O} }\textbf{u}'', {\bf{v}}''\right>_{\Omega}+  \left< {\bf{h}}_E {\bf{G}}_E \phi_E', \psi_E'\right>_{\Omega} =\left<U,V\right>_{\mc{X}_1}.~~~~\end{eqnarray*}
 Let $\mc{X}_{-1}$ denote the dual of $\mc{X}_1$ with respect to $\mc{X}$.
 By the Lax-Milgram theorem, $A_1$  extends to  an isomorphism  between $\mc{X}_1$ and $\mc{X}_{-1}.$    Therefore,  the inner product on $\mc{X}$ extends continuously to the duality pairing  $\left<\cdot,\cdot\right>_{\mc{X}_{-1}, \mc{X}_1}$ which satisfies (for $U,V\in \mc{X}_1$)
 $$\left<A_1 U, V\right>_{\mc{X}_{-1}, \mc{X}_1}= a(U'; V')= K  \left<u''', v'''\right>_{\Omega}+ \left< {\bf{h}}_{\mc O} {\bf{E}}_{\mc O} {\bf{u}}'',   {\bf{v}}''\right>_{\Omega}+   \left<{\bf{G}}_E {\bf{h}}_E \phi_E', \psi_E'\right>_{\Omega}$$
for the (h-N) boundary conditions and
$$\left<A_1 U, V\right>_{\mc{X}_{-1}, \mc{X}_1}= a(U; V)= K  \left<u'', v''\right>_{\Omega}+ \left< {\bf{h}}_{\mc O} {\bf{E}}_{\mc O} {\bf{u}}',   {\bf{v}}'\right>_{\Omega}+   \left<{\bf{G}}_E {\bf{h}}_E \phi_E, \psi_E\right>_{\Omega}$$
for the (c-D) and (m-m) boundary conditions.
Furthermore, we have dense compact embeddings $\mc{X}_1 \hookrightarrow \mc{X}\hookrightarrow \mc{X}_{-1}. $
\vspace{0.1in}

   From (\ref{eq2001}),  $A_1$ is a positive and self-adjoint operator.  Therefore there exists a sequence of orthogonal eigenvectors $\{E_{k,l}\}\in \mc{X}_1, k\ge 1, 1\le l\le m_k$  corresponding to the eigenvalues $\lambda_k$  and
    \begin{eqnarray}
\nonumber  && A_1E_{k,l}=\lambda_k E_{k,l}, ~~ 1\le l\le m_k \\
\label{eigen} \lambda_k>0,~~\lambda_k\to \infty, &&~ 1\le l\le m_k~~ {\rm{as}}~ k\to \infty,~~~ E_{k,l}\perp E_{m,n} ~{\rm{if}}~ k \ne m.
\end{eqnarray}
By (\ref{eq2001}), we have
\begin{eqnarray}
\nonumber \left<A_1 E_{k,l}, E_{k,l}\right>_{\mc{X}}=\left<\lambda_k E_{k,l}, E_{k,l}\right>_{\mc{X}}=\lambda_k \|E_{k,l}\|^2_{\mc{X}}=\|E_{k,l}\|^2_{\mc{X}_1}. \label{eq2002}\end{eqnarray}
Every $U\in \mc{X}_1$ has a unique orthogonal expansion $\sum_{k\ge 1, 1\le l\le m_k} {c_{k,l} E_{k,l}}$ and it follows from (\ref{eq2001}) that we have
\begin{eqnarray}
 \|U\|^2_{\mc{X}_1}=\sum_{k\ge 1, 1\le l\le m_k} \|c_{k,l}E_{k,l}\|^2_{\mc{X}_1}=\sum_{k\ge 1, 1\le l\le m_k} \lambda_k c^2_{k,l} \|E_{k,l}\|^2_{\mc{X}}.~~~~~~~~
\label{evestimate}
\end{eqnarray}

The inner product on $\mc{X}_{-1}$ is defined by
\begin{eqnarray}\label{innerINX-1}\left<U,V\right>_{\mc{X}_{-1}}=\left<A_1^{-1}U,A_1^{-1}V\right>_{\mc{X}_1}.
\end{eqnarray}
Note that the eigenfunctions $\{E_{k,l}\}_{ k\ge 1, 1\le l\le m_k}$ preserves their orthogonality in $\mc{X}$ and $\mc{X}_{-1}.$ Therefore, every $U\in \mc{X}~ ({\rm or}~\mc{X}_{-1})$ has a unique orthogonal expansion of the form $\sum_{k\ge 1, 1\le l\le m_k} {c_{k,l} E_{k,l}}$
   converging in $\mc{X}~ ({\rm or} ~~\mc{X}_{-1}),$ and we have
$$\|U\|^2_{\mc{X}}=\sum_{k\ge 1, 1\le l\le m_k} c^2_{k,l}\|E_{k,l}\|^2_{\mc{X}},$$
and respectively \begin{eqnarray}
\nonumber \|U\|^2_{\mc{X}_{-1}}&&=\sum_{k\ge 1, 1\le l\le m_k} c^2_{k,l}\|E_{k,l}\|^2_{\mc{X}_{-1}}=\sum_{k\ge 1, 1\le l\le m_k} c^2_{k,l}\|A_1^{-1}E_{k,l}\|^2_{\mc{X}_1}\\
\label{evestimate1} && =\sum_{k\ge 1, 1\le l\le m_k} \lambda_k^{-2} c^2_{k,l}\|E_{k,l}\|^2_{\mc{X}_1}=\sum_{k\ge 1, 1\le l\le m_k} \lambda_k^{-1} c^2_{k,l}\|E_{k,l}\|^2_{\mc{X}}.
\end{eqnarray}
 Eq. (\ref{evestimate1}) provides one  characterization of $\mc{X}_{-1}$.       However, we would like a function space characterization, particularly in the case of (m-m) boundary conditions.

We will need to refer Lemmata \ref{lem*} and \ref{dualch} below, which are proved in \cite{O-Hansen1}, and are adaptations of similar results in \cite{LT2}.
\vspace{0.1in}

\begin{Lem}\label{lem*} Let
$
 \mathrm{H} ={\rm{span}} \left\{ \sinh{\frac{x-L}{\sqrt{\alpha/m}}}\right\}\subset L^2(\Omega).
$ Let $\mc L$ be the operator $mI-\alpha D_x^2$ on the domain $H^2(\Omega)\cap H^1_0(\Omega).$ Then the restriction of $\mc L$ to $H^2_\#(\Omega)$  is an isomorphism from $H^2_\#(\Omega)$ to $\mathrm{H}^\perp$ in $L^2(\Omega).$
\end{Lem}
\vspace{0.1in}

 \begin{Lem} \label{dualch}$\mathrm{H}^\perp=(L^2(\Omega)/\mathrm{H})',$  where the duality is with respect to the $L^2(\Omega)$ inner product.
 \end{Lem}
\vspace{0.1in}

Now consider specifically the (m-m) boundary conditions.   For $V= (v,\mathbf{v}) \in \mc{X}_1= H^2_\#(\Omega) \times ( H^1_\dag(\Omega))^{(m+1)} ,\ U= (u, \mathbf{u}) \in \mc{X}=H^1_0(\Omega)\times (L^2(\Omega))^{(m+1)} $,  an integration by parts of
(\ref{eq2003}) results in
\begin{equation}  \nonumber c(U,V) = -\left<u, \mc L v\right>_{\Omega}+ \left< {\bf{h}}_{\mc O} {\bf{p}}_{\mc O} {\bf{u}},   {\bf{v}}\right>_{\Omega}.\end{equation}
The second term remains bounded for all ${\bf{u}}\in (H^1_\dag(\Omega))^{(m+1)})' $  (with duality relative to $L^2(\Omega)$).      In the first term, however,
by Lemma \ref{lem*},  the range of $\mc L$  is $\mathrm{H}^\perp$ in $L^2(\Omega).$   Hence for the first term to remain bounded,  by Lemma \ref{dualch},
$u\in L^2(\Omega)/\mathrm{H}$.   Therefore, in the case of (m-m) boundary conditions,
\begin{equation}  \mc{X}_{-1} =  L^2(\Omega)/\mathrm{H}\times (H^1_\dag(\Omega))^{(m+1)})'  \label{xm1}\end{equation}

 It is easiest to
characterize  $\mc{H}_{-1}$ in the undamped case.   (Later we will show that the same characterization holds in the damped case.)
Write the operator $\mc{A}$ as follows:
\begin{equation} \nonumber \mc{A} = \mc{A}_0 +  \mc{B} =  \left(\begin{array}{cc} 0&I \\ -A_1 & 0\end{array}\right)  + \left(\begin{array}{cc} 0& 0 \\0 & A_2\end{array} \right)
\end{equation}
Then $\mc{D}(\mc{A}) =\mc{D}(\mc{A}_0)$ and hence
$\mc{A}_0: \mc{H}= \mc{X}_1\times \mc{X} \to\mc{H}_{-1}$ is an isomorphism by Theorem \ref{eigens} and Corollary \ref{ext}.
It follows that an  inner product on $\mc{H}_{-1}$ can be defined by  $\left<Y,Z\right>_{\mc{H}_{-1}}=\left<  { \mc {A}_0}^{-1}Y,  { \mc{A}_0}^{-1} Z\right>_{\mc H}$.   Hence,   in the undamped case,
\begin{eqnarray}\nonumber \left<Y,Z\right>_{\mc{H}_{-1}}&=&\left< { \mc {A}_0}^{-1}Y, { \mc{A}_0}^{-1} Z\right>_{\mc H}\\
\nonumber &=& c(Y_1, Z_1)+ a( -A_1^{-1}Y_2,  -A_1^{-1}Z_2)\\
\nonumber &=& \left<Y_1, Z_1\right>_{\mc X}+ \left< A_1^{-1} Y_2,  A_1^{-1} Z_2\right>_{\mc X_1} \\
\nonumber &=& \left<Y_1, Z_1\right>_{\mc X}+ \left< Y_2, Z_2\right>_{\mc X_{-1}}
\end{eqnarray}
where we used (\ref{inner}) and (\ref{innerINX-1}). By (\ref{xm1}),  we have  in the undamped case with (m-m) boundary conditions,
\begin{eqnarray}\label{H-1}\mc{H}_{-1}=\mc X\times \mc X_{-1}=H^1_0(\Omega)\times (L^2(\Omega))^{(m+1)}\times (L^2(\Omega)/\mathrm{H}) \times (H^1_\dag(\Omega)')^{(m+1)}.\quad\quad
\end{eqnarray}

\section{Observability results and the Proof of Theorem \ref{observability}} \label{C2S3-Observ}
\vspace{0.1in}
We prove our main observability results in this section.     We begin  with some preliminary results for the decoupled system.
\vspace{0.1in}

\subsection{{Observability results for  decoupled system}}
\label{sub1}

 Consider (\ref{mainhomo}) without the coupling terms, i.e.,  with ${\bf{G}}_E=\tilde {\bf{G}}_E = 0.$    What remains is a Rayleigh beam equation
 and $(m+1)$ wave equations:
\begin{equation}\left\{ \begin{array}{l}
m\ddot z -\alpha \ddot z'' +  K   z''''  = 0 ~~~ {\rm{on}}~~ \Omega\times \mathbb{R}^+ \\
     {\ddot v}_{\mc O} -{\bf{p}}_{\mc O}^{-1} {\bf{E}}_{\mc O}  {v}_{\mc O}'' = 0~~ {\rm{on}} ~~ \Omega\times \mathbb{R}^+,
 \end{array} \right.
\label{mainhomouc}
\end{equation}
 with the boundary conditions (\ref{bccc}) and the initial conditions (\ref{initialhomo}).
 Let $$U=:(u,{\bf{u}})=(z, {v}_{\mc O})^{\rm T}, ~~~~V:= (v, {\bf{v}})^{\rm T}=(\dot z,  {\dot v}_{\mc O})^{\rm T}, ~~{\rm{and}} ~~ Y:=(U,V)^{\rm T}.$$  Then the semigroup corresponding to (\ref{mainhomouc}) is given by
\begin{eqnarray}\nonumber &&\frac{dY}{dt}=\mc{A}_d Y:=\left( {\begin{array}{*{20}c}
   {0} & {I }  \\
   {-A_d } & {0} &   \\
\end{array}} \right) \left( \begin{array}{l}
 U \\
 V \\
 \end{array} \right), \\
\nonumber  && Y(0)=(U(0),V(0))^{\rm T}=(z^0,  {v}^0_{\mc O}, z^1,  {v}^1_{\mc O} )^{\rm T}
\end{eqnarray}
 where $A_d U:=\left( {\begin{array}{*{20}c}
   K \mc L^{-1} u''''  \\
   {  {-\bf{p}}^{-1}_{\mc O} } {  {\bf{E}}_{\mc O} }\textbf{u}'' \\
\end{array}} \right)\label{A1d}.$ Define the quadratic forms $a_d$ and $c_d$ by
\begin{eqnarray}
\nonumber && c_d(z, v_{\mc O} )= m\left<z,z\right>_{\Omega}+\alpha\left<z',z'\right>_{\Omega}+ \left< {\bf{h}}_{\mc O} {\bf{p}}_{\mc O} { v}_{\mc O},   { v}_{\mc O}\right>_{\Omega}\\
\nonumber && a_d(z, v_{\mc O} )= K  \left<z'',z''\right>_{\Omega}+ \left< {\bf{h}}_{\mc O} {\bf{E}}_{\mc O} { v}_{\mc O}',   { v}_{\mc O}'\right>_{\Omega}.
\end{eqnarray}
 The natural and ``higher order'' energies of the decoupled system are given by
 \begin{subequations}
\begin{empheq}[left={\phantomword[r]{}{\mc{E}_d(t) = }  \empheqlbrace}]{align}
 \nonumber &\frac{1}{2}\left( a_d(z', v_{\mc O}')+c_d(\dot z', \dot v_{\mc O}')\right) & \text{(h-N)}\\
\nonumber & \frac{1}{2}\left(a_d( z, v_{\mc O})+c_d(\dot z, \dot v_{\mc O})\right). &  \text{(c-D), (m-m)}.
\end{empheq}
\end{subequations}
The energy inner products corresponding to each set of boundary conditions are defined by
 \begin{subequations}
\begin{empheq}[left={\phantomword[r]{}{\left<Y,\widehat{Y}\right>_{\mc{H}}=}  \empheqlbrace}]{align}
\nonumber & ~a_d(U';\widehat{U}')+c_d(V';\widehat{V}'). & \text{(h-N)}\\
\nonumber & ~a_d(U;\widehat{U})+c_d(V;\widehat{V}) & \text{(c-D)}, \text{(m-m)}.
\end{empheq}
\end{subequations}

 In the above $\mc{A}_d$ is densely defined by $\mc{A}_d: \mc{D}(\mc{A}_d) \subset \mc{H}\to \mc{H} $  and  note that
$\mc{D}(\mc{A}_d)=\mc{D}(\mc{A}).$
\vspace{0.1in}

   \begin{rmk} \label{equivalence} (i) It is easy to verify that $\mc{E}(t)\asymp \mc{E}_d(t), ~\forall t>0.$ Indeed, for the hinged-Neumann (h-N) boundary conditions
   \begin{eqnarray}\nonumber
   \nonumber \left|\left<{\bf{G}}_E {\bf{h}}_E \phi_E',\phi_E'\right>_{\Omega}\right| = && \left|\left<{\bf{G}}_E {\bf{h}}_E^{-1}\left({\bf{B}} { v}_{\mc O}' + {\bf{h}}_E  N z''\right),( {\bf{B}} { v}_{\mc O}' + {\bf{h}}_EN z''\right>_{\Omega}\right|\\
   \nonumber \le && ~C\left(\|v_{\mc O}''\|^2_{(L^2(\Omega))^{(m+1)}} + \|z'''\|^2_{L^2(\Omega)}\right)\le  C\mc{E}_d,
   \end{eqnarray}
   and for the clamped-Dirichlet (c-D) and  mixed-mixed (m-m) boundary conditions
   \begin{eqnarray}\nonumber
   \nonumber \left|\left<{\bf{G}}_E {\bf{h}}_E \phi_E,\phi_E\right>_{\Omega}\right| = && \left|\left<{\bf{G}}_E {\bf{h}}_E^{-1}\left({\bf{B}} { v}_{\mc O} + {\bf{h}}_E  N z'\right),( {\bf{B}} { v}_{\mc O} + {\bf{h}}_EN z'\right>_{\Omega}\right|\\
   \nonumber \le && ~C\left(\|v_{\mc O}'\|^2_{(L^2(\Omega))^{(m+1)}} + \|z''\|^2_{L^2(\Omega)}\right)\le  C\mc{E}_d
   \end{eqnarray}  where $ C$ denotes a generic constant. Therefore,
   \begin{eqnarray}
   \label{energy2} \mc{E}_d\le \mc{E} \le C \mc{E}_d.
   \end{eqnarray}
   (ii) In the case of (m-m) boundary conditions, we define the solutions of (\ref{mainhomouc}),(\ref{bccc}) and (\ref{initialhomo})  on the extended space $\mc H_{-1}$ (defined by (\ref{H-1}))   in exactly the same way as we did for the undamped coupled system,  i.e., by applying  Corollary \ref{ext},  Lemma \ref{lem*}, and Lemma \ref{dualch} to the decoupled system.              Therefore we define the energy of the weak  solutions  by
   \begin{eqnarray}
\nonumber \mc{E}_{-1}(t)&=&  \frac{1}{2}\|(z, \dot z, v_{\mc O},\dot v_{\mc O})\|^2_{\mc{H}_{-1}}\\
\label{ext-energy}&\approx&  \frac{1}{2}\left(\|z\|^2_{H^1_0(\Omega)}+ \|\dot z\|^2_{(L^2(\Omega))^{(m+1)}}+ \|v_{\mc O}\|^2_{L^2(\Omega)/\mathrm{H}} + \|{\dot v_{\mc O}}\|^2_{((H^1_\dag(\Omega))')^{(m+1)}} \right).\quad\quad\quad
\end{eqnarray}

   \end{rmk}
   \vspace{0.1in}

The following results for the interior regularity, hidden regularity, and observability of the decoupled system (\ref{mainhomouc}) follow from the standard semigroup theory, standard results for the wave equation, e.g. see \cite{Komornik-P}, \cite{LT4}, and observability results obtained in \cite{O-Hansen1}.
\vspace{0.1in}

 \begin{thm} \label{thmdc}

 \noindent \textbf{\text{(a)}}  Consider
 \begin{equation} \left\{ \begin{array}{l}
m\ddot z-\alpha  \ddot z''+K  z'''' +f(x,t) = 0 \quad  {\rm{in}} ~~\Omega \times \mathbb{R}^+\\
{\ddot v}_{\mc O} -{\bf{p}}_{\mc O}^{-1} {\bf{E}}_{\mc O}  {v}_{\mc O}'' +f_{\mc O}(x,t) = 0 \quad {\rm{in}} ~~ \Omega\times \mathbb{R}^+
\label{raohomodecf}
   \end{array} \right.
\end{equation}
with the boundary conditions (\ref{bccc}) and the initial conditions
   \begin{equation} \nonumber z(x,0)=\dot z(x,0)=0, ~{v}_{\mc O}(x,0)= {\dot v}_{\mc O}(x,0)= 0 ~{\mbox{on}} ~\Omega. \end{equation}
Assume
\begin{subequations}
\begin{empheq}[left={\phantomword[l]{}{ }  \empheqlbrace}]{align}
\nonumber & f\in L^1(0,T; L^2(\Omega)), ~f_{\mc O}\in L^1(0,T; (H^1(\Omega))^{(m+1)}) & \text{(h-N)}\\
\nonumber & f\in L^1(0,T; H^{-1}(\Omega)),~f_{\mc O}\in L^1(0,T; (L^2(\Omega))^{(m+1)})   & \text{(c-D)}\\
\nonumber & f\in L^1(0,T; L^2(\Omega)/\mathrm{H})),~f_{\mc O}\in L^1(0,T; ((H^1_\dag(\Omega))')^{(m+1)}) & \text{(m-m)}.
\end{empheq}
\end{subequations}
Then
$(z, \dot z, v_{\mc O}, \dot v_{\mc O})\in C\left([0,T]; \mc{H}\right)$ and the solution of (\ref{raohomodecf})
satisfy for every $T>0$ the direct inequality
\begin{eqnarray}
\nonumber &&\int_{0}^T \left(|z'''(L,t)|^2+|v_{\mc O}''(L,t)|^2\right) ~dt  \le  C \|(f, f_{\mc O}')\|^2_{L^1\left(0,T;L^2(\Omega)\times (L^2(\Omega))^{(m+1)}\right)} \quad\\
\nonumber &&\int_{0}^T \left(|z''(L,t)|^2 + |v_{\mc O}'(L,t)|^2\right) ~dt \le C \|(f, f_{\mc O})\|^2_{L^1\left(0,T;H^{-1}(\Omega)\times (L^2(\Omega))^{(m+1)}\right)} \quad\\
\nonumber &&\int_{0}^T \left(|z'(L,t)|^2 + |v_{\mc O}(L,t)|^2 \right) ~dt \le C \|(f, f_{\mc O})\|^2_{L^1\left(0,T;L^2(\Omega)/\mathrm{H})\times ((H^1_\dag(\Omega))')^{(m+1)}\right)}\quad\quad\quad
\end{eqnarray}
 for (h-N), (c-D), and (m-m) respectively. In the above $C=C(T)$ is a generic constant.
\vspace{0.1in}

\noindent {\text{\textbf{(b)}}} Consider  \begin{equation} \left\{ \begin{array}{l}
m\ddot z-\alpha  \ddot z''+K  z''''  = 0  \quad {\rm{in}} ~~~~\Omega \times \mathbb{R}^+\\
{\ddot v}_{\mc O} -{\bf{p}}_{\mc O}^{-1} {\bf{E}}_{\mc O}  {v}_{\mc O}''  = 0\quad {\rm{in}} ~~ \Omega\times \mathbb{R}^+
\label{raohomodec}
   \end{array} \right.
\end{equation}
with the boundary conditions (\ref{bccc}) and the initial conditions  (\ref{initialhomo}). Assume that the initial conditions satisfy $(z_0,z_1, v_{\mc O}^0, v_{\mc O}^1)\in \mc{H}.$ Then $(z, \dot z, v_{\mc O},\dot v_{\mc O})\in C\left([0,T]; \mc{H}\right)$ and the solution of (\ref{raohomodec})
satisfies for every $T>\tau$ ($\tau$ is defined by (\ref{tau})) the following observability and hidden regularity results
\begin{eqnarray}
\nonumber \int_{0}^T \left(|z'''(L,t)|^2+|v_{\mc O}''(L,t)|^2 \right)~dt  \asymp  \mc{E}_d(0) && ~~~\text{(h-N)}\\
\nonumber \int_{0}^T \left( |z''(L,t)|^2 + |v_{\mc O}'(L,t)|^2 \right)~dt \asymp  \mc{E}_d(0) &&~~~ \text{(c-D)}\\
\nonumber \int_{0}^T \left(|z'(L,t)|^2 + |v_{\mc O}(L,t)|^2 \right)~dt  \asymp  \mc{E}_{-1}(0) && ~~~\text{(m-m)}
\end{eqnarray}
where $\mc E_{-1}$ is defined by (\ref{ext-energy}).
 \end{thm}

\vspace{0.1in}
\subsection{Observability results for  coupled, undamped system}         We now consider the coupled, undamped system , i.e. ${\bf{G}}_E \ne 0, ~\tilde {\bf{G}}_E = 0$.
Consider (\ref{mainhomo}) without the damping terms, i.e., $\tilde {\bf{G}}_E = 0$:
\begin{equation}\left\{ \begin{array}{l}
m\ddot z -\alpha \ddot z'' +  K   z'''' -   N^{\rm T} {\bf{h}}_E {\bf{ G}}_E \phi_E' = 0 ~~~ {\rm{on}}~~ \Omega\times \mathbb{R}^+ \\
    {\ddot v}_{\mc O} -{\bf{p}}_{\mc O}^{-1} {\bf{E}}_{\mc O}  {v}_{\mc O}'' +  {\bf{p}}_{\mc O}^{-1}{\bf{h}}_{\mc O}^{-1}{\bf{B}}^{\rm T} {\bf{ G}}_E \phi_E  = 0~~ {\rm{on}} ~~ \Omega\times \mathbb{R}^+ \\
 \text{where}~( {\bf{B}} { v}_{\mc O}={\bf{h}}_E \phi_E-{\bf{h}}_E  N z')
 \end{array} \right.
\label{mainhomoud}
\end{equation}
with the boundary conditions (\ref{bccc}) and the initial conditions (\ref{initialhomo}).    Since  the generator $\mc A_0$ is skew-adjoint,  the energy
  $\mc E$  in  (\ref{energy})  is conserved along solution trajectories.

Now consider the eigenvalue problem corresponding to (\ref{mainhomoud})

\begin{equation}\label{acik}\mc{A}_0
\left( \begin{array}{l}
 U \\
 V \\
 \end{array} \right) = \lambda \left( \begin{array}{l}
 U \\
 V \\
 \end{array} \right) \Rightarrow ~~ V=\lambda U ~~~~and ~~~~  A_1 U =\lambda V.
 \end{equation}
Explicitly, (\ref{acik}) can be written as
 \begin{subequations}\label{ep}
\begin{empheq}[left={\phantomword[l]{}{ } \empheqlbrace}~~]{align}
\label{ep1}& -Ku''''+ N^{\rm T} {\bf{h}}_E {\bf{G}}_E\phi_E '= \lambda^2 \mc L u   & \\
\label{ep2}& { {\bf{h}}_{\mc O} {\bf{E}}_{\mc O} }{\textbf{u}}''-{\bf{B}}^{\rm T} {\bf{G}}_E  \phi_E  = \lambda^2  { {\bf{p}}_{\mc O} {\bf{h}}_{\mc O}  }{\textbf{u}}. &
\end{empheq}
\end{subequations}
The following is the {\it key uniqueness result} of this paper.

\vspace{0.1in}
\begin{Lem} \label{xyz} The eigenvalue problem \label{overtrivial}(\ref{ep}) together with  any of the following sets of boundary conditions
 the boundary conditions
\begin{eqnarray}
 &  \left\{ \begin{array}{l}
  u(0,t)=u''(0,t)=u(L,t)=u''(L,t)=u'''(L,t)=0 \label{bdry3udov}\\
 { \bf u}'(0,t)={ \bf u}'(L,t)={ \bf u}''(L,t)=0, ~~~\label{bdry4udov}
   \end{array} \right\} & ~~~\text{(h-N)}\\
 &  \left\{ \begin{array}{l}
  u(0,t)=u'(0,t)=u(L,t)=u'(L,t)=u''(L,t)=0 \label{bdry5udov}\\
 { \bf u}(0,t)={ \bf u}(L,t)={ \bf u}'(L,t)=0, ~~~ \label{bdry6udov}
   \end{array} \right\} & ~~~\text{(c-D)},  \text{(m-m)}\quad\quad
\end{eqnarray}
has only  the trivial solution.
\end{Lem}
\vspace{0.1in}

 \textbf{Proof:} We first consider the case of (h-N) boundary conditions. Note that if $(u,\bf u)$ satisfies (\ref{ep})-(\ref{bdry4udov}), then $(z,{\bf z})=(u'',{\bf u}'')$ satisfies (\ref{ep}) with the boundary conditions
  \begin{subequations}\label{overlamd}
\begin{empheq}[left={\phantomword[l]{}{ } \empheqlbrace}~~]{align}
\label{overlamd1}&  z(0,t)=z''(0,t)= z(L,t)=z'(L,t)=z''(L,t) = 0    & \\
\label{overlamd2}&  \textbf{z}'(0,t) = \textbf{z}(L,t) = \textbf{z}'(L,t)=  0.  &
\end{empheq}
\end{subequations}
 If $(z,{\bf z})\equiv 0,$ then  $(u'',{\bf u}'')\equiv 0$ by using the boundary conditions (\ref{bdry4udov}).   Thus in any of the cases, it is enough to show that (\ref{ep}),(\ref{bdry4udov}) and (\ref{ep}),(\ref{bdry6udov}) have only the trivial solutions.

Now multiply (\ref{ep1}) by $x\bar u'-3 \bar u$ and multiply (dot product) (\ref{ep2}) by $x\bar{\textbf{u}}'-2\bar{\textbf{u}}$ respectively and add to each other. Then integrating by parts on $\Omega$ with the use of boundary conditions (\ref{overlamd})  yields :
\begin{eqnarray}
\nonumber  && 0=\int_{\Omega} \left(  -4\lambda^2|u|^2 - 2\alpha\lambda^2 |u'|^2    -3 \lambda^2  {\bf{h}}_{\mc O}   {\bf{p}}_{\mc O}  \textbf{u}\cdot   \bar{\textbf{u}}-{ {\bf{h}}_{\mc O} {\bf{E}}_{\mc O} } {\textbf{u}}' \cdot  \bar{\textbf{u}}' \right)~ dx\\
\nonumber && ~~~~ +\int_{\Omega} \left( -x \lambda^2 \bar u u'  +\alpha x\lambda^2 u' \bar u''   -K x \bar u''''xu'  -\lambda^2 {\bf{h}}_{\mc O}   {\bf{p}}_{\mc O}  \textbf{u}' \cdot  x\bar{\textbf{u}}\right) ~dx \\
 \label{multres11}&&  ~~~~+\int_{\Omega} \left({ {\bf{h}}_{\mc O} {\bf{E}}_{\mc O} } {\textbf{u}}' \cdot    x\bar{\textbf{u}}''-3{\bf{G}}_E \phi_E \cdot h_E\bar \phi_E - {\bf{G}}_E \phi_E' \cdot x h_E\bar \phi_E\right)   ~dx.
 \end{eqnarray}
Now we look at the solution $(\bar u, \bar{\textbf{u}})$ of the eigenvalue problem (\ref{ep}) corresponding to the eigenvalue $\bar \lambda:$
\begin{subequations}
\begin{empheq}[left={\phantomword[l]{}{ } \empheqlbrace}~~]{align}
\label{epcon1} &\bar\lambda^2  \bar u -\alpha \bar\lambda^2 \bar u'' + K\bar u'''' - N^{\rm T} {\bf{h}}_E {\bf{G}}_E \bar{\phi}_E' = 0 \\
\label{epcon2}&\bar \lambda^2  {\bf{h}}_{\mc O}   {\bf{p}}_{\mc O}  \bar{\textbf{u}} - { {\bf{h}}_{\mc O} {\bf{E}}_{\mc O} } \bar{{\textbf{u}}}''+{\bf{B}}^{\rm T} {\bf{G}}_E \bar{\phi}_E =0.
 \end{empheq}
\end{subequations} with the conjugate boundary conditions
  \begin{subequations}\label{overlamdcon}
\begin{empheq}[left={\phantomword[l]{}{ } \empheqlbrace}~~]{align}
\label{overlamdcon1}&  \bar u(0,t)=\bar u''(0,t)= \bar u(L,t)=\bar u'(L,t)=\bar u''(L,t) = 0    & \\
\label{overlamdcon2}& \bar {\textbf{u}}'(0,t) = \bar {\textbf{u}}(L,t) = \bar{\textbf{u}}'(L,t)=  0  &
\end{empheq}
\end{subequations}
 Now multiply (\ref{epcon1}) by $xu'+2u$ and multiply (dot product) (\ref{epcon2}) by $x\textbf{u}'+ 3 \textbf{u}$ respectively and add to each other. Then integrating by parts on $\Omega$ with the use of (\ref{overlamdcon}) yields
   \begin{eqnarray}
 \nonumber &&0= \int_{\Omega} \left( \bar\lambda^2  \bar u x u' -\alpha \bar\lambda^2 \bar u'' x u' + K \bar u'''' x  u' +\bar\lambda^2  {\bf{h}}_{\mc O}   {\bf{p}}_{\mc O}  \bar{\textbf{u}}\cdot  x {\textbf{u}}' - { {\bf{h}}_{\mc O} {\bf{E}}_{\mc O} } \bar{{\textbf{u}}}'' \cdot  x {\textbf{u}}'\right)~ dx\\
 \nonumber && \quad + \int_{\Omega} \left( 2\bar\lambda^2  |u|^2 +2\alpha \bar\lambda^2 |u'|^2 + 2K|u''|^2 +3\bar\lambda^2  {\bf{h}}_{\mc O}   {\bf{p}}_{\mc O}  \bar{\textbf{u}}\cdot   \textbf{u}+3{ {\bf{h}}_{\mc O} {\bf{E}}_{\mc O} } \bar{{\textbf{u}}}' \cdot  \textbf{u}' \right)~ dx\\
 \label{mult7-1} &&\quad + \int_{\Omega} \left(3{\bf{G}}_E \bar\phi_E\cdot h_E \phi_E dx + {\bf{G}}_E \bar \phi_E\cdot(x h_E \phi_E' )\right) ~dx.
  \end{eqnarray}
Eventually, adding (\ref{multres11}) and (\ref{mult7-1}) gives
\begin{eqnarray}
\nonumber && 0= \int_{\Omega} \left(  -2(2\lambda^2-\bar\lambda^2)|u|^2 -2\alpha(\lambda^2-\bar\lambda^2)|u'|^2 + 2K|u''|^2 \right) ~dx\\
\nonumber && \quad + \int_{\Omega}    \left(-3(\lambda^2-\bar\lambda^2){\bf{h}}_{\mc O}   {\bf{p}}_{\mc O}  \bar{\textbf{u}}\cdot {\textbf{u}}  + 2{ {\bf{h}}_{\mc O} {\bf{E}}_{\mc O} } {\textbf{u}}' \cdot  \bar{\textbf{u}}'\right)~ dx\\
\nonumber && \quad + \int_{\Omega}  \left( x (-\lambda^2+\bar\lambda^2) \bar u u'  +\alpha x(\lambda^2-\bar\lambda^2) u' \bar u''   +(-\lambda^2+\bar\lambda^2) {\bf{h}}_{\mc O}   {\bf{p}}_{\mc O}  \textbf{u}' \cdot  x\bar{\textbf{u}} \right)~ dx\\
 \label{multreslast}&& \quad  +\int_{\Omega}  \left( ({\bf{G}}_E \bar \phi_E )\cdot(x h_E \phi_E' ) - ({\bf{G}}_E  \phi_E )\cdot(x h_E \bar\phi_E' ) \right)~ dx.
 \end{eqnarray}
Note that energy of the undamped system is conserved. Therefore, all eigenvalues are located on the imaginary axis. Now let $\lambda=\mp is, ~s\in \mathbb{R}^+.$ Then  $\lambda^2$ and $\bar\lambda^2$ have the same sign.  Then (\ref{multreslast}) reduces to
\begin{eqnarray}
\nonumber && \int_{\Omega}  2s^2|u|^2 +  2K|u''|^2 + 2{ {\bf{h}}_{\mc O} {\bf{E}}_{\mc O} } {\textbf{u}}' \cdot  \bar{\textbf{u}}' ~dx\\
\label{multreslastt}&&  \quad +\int_{\Omega}  \left( ({\bf{G}}_E \bar \phi_E )\cdot(x h_E \phi_E' ) - ({\bf{G}}_E  \phi_E )\cdot(x h_E \bar\phi_E' ) \right)~ dx=0.
 \end{eqnarray}
    Note  that the last two terms in (\ref{multreslastt}) are conjugates of each other. Therefore the second integral term is pure imaginary. Hence we have $u''=0$ and $\textbf{u}'=0.$ Using boundary conditions (\ref{overlamd})  we get $(u, {\bf u})\equiv 0.$ This completes the proof for the (h-N) boundary conditions.

     In (c-D) and (m-m) cases, similar calculations again lead to (\ref{multreslastt}). Hence using boundary conditions (\ref{bdry6udov}), we obtain $(u, {\bf u})\equiv 0.$ $\square$
     \vspace{0.1in}

 The following result is  Theorem 6.2   in (Chap VI, \cite{Komornik-P}), as it applies to our problem.    \vspace{0.1in}

\begin{thm} \label{Haraux-type-thm}  Let $Y=[z, v_{\mc O}, \dot z, \dot v_{\mc O}]^{\rm T}$ and $Y_0=[z^0,v_{\mc O}^0, z^1, v_{\mc O}^1]^{\rm T}.$  Assume
the following two conditions.\vspace{0.1in}

  \noindent (i) There exists a sufficiently large $k'\in \mathbb{N}$ such that for $T>\tau$ ($\tau$ is defined by (\ref{tau})) we have
\begin{subequations}\label{partI}
\begin{empheq}[left={\phantomword[l]{}{ } \empheqlbrace}]{align}
  \label{parti} {\int_0^T \left(|z'''(L,t)|^2+| v_{\mc O}''(L,t)|^2\right) ~dt}\asymp \|Y_0\|_{\mc H}^2&&~~~ \text{(h-N)}\\
  \label{parti2}\int_{0}^T \left( |z''(L,t)|^2 + |v_{\mc O}'(L,t)|^2 \right)~dt \asymp \|Y_0\|_{\mc H}^2 &&~~~ \text{(c-D)}\\
\label{parti3}\int_{0}^T \left(|z'(L,t)|^2 + |v_{\mc O}(L,t)|^2 \right)~dt  \asymp \|Y_0\|_{\mc H_{-1}}^2 && ~~~\text{(m-m)}
\end{empheq}
\end{subequations}
for all solutions of (\ref{mainhomoud}) with $Y_0\in _{\mc{H}_{k'}^\perp}$ where $\mc{H}_{k'}= {\rm{span}} \{E_{k,l}, ~1\le k\le k',~ 1\le l\le m_k\}.$
\vspace{0.1in}

  \noindent (ii) There exists $\bar T >0$ such that for all $T>\bar T$ the estimates (\ref{partI})
  hold for all solutions of (\ref{mainhomoud}) with $Y_0$ such that $\mc{A}Y_0=\lambda Y_0.$
\vspace{0.1in}

  Then for any $T>\tau $ the estimates (\ref{partI})  hold for all solutions $Y_0\in \mc{H}$ for the (h-N) and (c-D) cases, and $Y_0\in \mc H_{-1}$ for the (m-m) case.
\end{thm}
\vspace{0.1in}

We are now able to prove our main observability result  (Theorem \ref{observability}) for the undamped system (with $\tilde {\bf{G}}_E\equiv 0$):
\vspace{0.1in}

\begin{Lem}  \label{aaa}  Let $T>\tau$, where $\tau$ is given by (\ref{tau}) and assume that  $\tilde {\bf{G}}_E\equiv 0.$   Then solutions of  (\ref{mainhomoud})  satisfy the observability and hidden regularity estimates  (\ref{obs}).     \end{Lem}
\vspace{0.1in}

\textbf{Proof:}   This will follow from Theorem \ref{Haraux-type-thm} once we verify the conditions (i) and (ii) of the hypothesis are satisfied.

 First we consider the case of (h-N) boundary conditions. Let us write the solution of (\ref{mainhomoud}) in the form
 $$(z,v_{\mc O})^{\rm T}=(z_f,{v_{\mc O}}_f)^{\rm T} +(\hat{z}, \hat{v}_{\mc O})^{\rm T}.$$
 where $[z_f,{v_{\mc O}}_f]^{\rm T}$ solves (\ref{raohomodecf}) with
 $$(f,f_{\mc O})^{\rm T}=[-N^{\rm T} {\bf{h}}_E {\bf{ G}}_E \phi_E',{\bf{p}}_{\mc O}^{-1}{\bf{h}}_{\mc O}^{-1}{\bf{B}}^{\rm T}  {\bf{ G}}_E \phi_E]^{\rm T}$$ and zero initial conditions, and $(\hat{z}, \hat{v}_{\mc O})^{\rm T}$ solves (\ref{raohomodec}) with the initial data $ (z^0, v_{\mc O}^0, z^1, v_{\mc O}^1)^{\rm T}$ where $ {\bf{B}} { v}_{\mc O}={\bf{h}}_E \phi_E-{\bf{h}}_E  N z'.$  For $T>\tau,$ we apply  part (a) of Theorem \ref{thmdc} for $(z_f,v_{{O}_f})^{\rm T},$ and obtain
    \begin{eqnarray}
\nonumber &&  \int_0^T \left( |  z_f'''(L,t)|^2 + |v_{{\mc O}_f}''(L,t)|^2 \right)~ dt  \\
\nonumber && \quad\quad\quad\le \int_0^T \left(\| N^{\rm T} {\bf{h}}_E{\bf{ G}}_E {\bf{B}}   v_{\mc O}'\|^2_{L^2(\Omega)} +\| N^{\rm T} {\bf{h}}_E {\bf{ G}}_E {\bf{h}}_E N   z''\|^2_{L^2(\Omega)}\right.  \\
\nonumber &&  \quad\quad\quad\quad +  \left.\|{\bf{p}}_{\mc O}^{-1}{\bf{h}}_{\mc O}^{-1}{\bf{B}}^{\rm T} {\bf{ G}}_E {\bf{h}}^{-1}_E {\bf{B}} v_{\mc O}'  \|^2_{(L^2(\Omega))^{m+1}} +  \|{\bf{p}}_{\mc O}^{-1}{\bf{h}}_{\mc O}^{-1}{\bf{B}}^{\rm T} {\bf{ G}}_E N  z''  \|^2_{(L^2(\Omega))^{m+1}}\right)~ dt
\end{eqnarray}
and therefore
\begin{eqnarray}
\nonumber  &&  \int_0^T  \left(|  z_f'''(L,t)|^2 + |v_{{\mc O}_f}''(L,t)|^2\right)~  dt  \le  C_1({\bf{ G}}_E) \int_0^T \left(\|   v_{\mc O}'\|^2_{(L^2(\Omega))^{m+1}} +  \|   z''\|^2_{L^2(\Omega)} \right)~ dt
\end{eqnarray}
where $C_1$ is a function of ${\bf{ G}}_E.$  It follows from (\ref{evestimate}) that
\begin{eqnarray}\label{ineq}\|(z, v_{\mc O})^{\rm T}\|_{\mc{X}_1}^2\ge \lambda_1 \|(z, v_{\mc O})^{\rm T}\|_{\mc{X}}^2,\end{eqnarray}
where  $\{\lambda_k\}_{k=1 }^{\infty}$ are the eigenvalues of the operator $A_1.$   By equivalence of the energy (see Remark \ref{equivalence}) and (\ref{ineq})
 it follows that
\begin{eqnarray}
\nonumber  &&\int_0^T  \left(|  z_f'''(L,t)|^2 + |v_{{\mc O}_f}''(L,t)|^2\right)~ dt\\
\label{eq1001} && \quad \le   C_2({\bf{ G}}_E) \int_0^T \left(\frac{1}{\sqrt{\lambda_{1}}}\|    v_{\mc O}''\|^2_{(L^2(\Omega))^{m+1}} + \frac{1}{\sqrt{\lambda_{1}}}\|   z'''\|^2_{L^2(\Omega)}\right)~ dt \le   \frac{C_3({\bf{ G}}_E)}{\sqrt{\lambda_{1}}}  \mc{E}_d(0).\quad\quad\quad
\end{eqnarray}
Now if we use the assumption $Y_0 \perp \{E_{k,l}, ~1\le k\le k',~ 1\le l\le m_k\},$  in part (i) of the theorem, then
we have
\begin{equation}\|(z, v_{\mc O})^{\rm T}\|^2_{\mc{X}_1}\ge \lambda_k' \|(z, v_{\mc O})^{\rm T}\|^2_{\mc{X}}\label{important-ineq}
\end{equation}
and therefore (\ref{eq1001}) can be written as
  \begin{eqnarray} \label{eq1002}  \int_0^T  |  z_f'''(L,t)|^2 + |v_{{\mc O}_f}''(L,t)|^2 dt   \le   \frac{C_3({\bf{ G}}_E)}{\sqrt{\lambda_{k'}}}   \mc{E}_d(0).
\end{eqnarray}

 Next, for $T>\tau$  we apply  part (b) of Theorem  \ref{thmdc} together with (\ref{raohomodec})  for $(\hat { z},\hat{y}_{\mc O})^{\rm T}$  respectively, for $c_1,c_2>0$ we get
 \begin{eqnarray}
\label{eq1003} && c_1 \mc{E}_d(0)\le \int_0^T  |  \hat{z}'''(L,t)|^2 + |\hat{v}_{\mc O}''(L,t)|^2 dt  \le c_2  \mc{E}_d(0).
\end{eqnarray}
Since
\begin{eqnarray} \label{eq1015}  |z'''|^2  \le 2 | \hat{z}'''|^2 + 2| z_f'''|^2, \quad
| v_{\mc O}''|^2 \le 2 | \hat{v}_{\mc O}''|^2 + 2| v_{{\mc O}_f}''|^2
\end{eqnarray}
By combining (\ref{eq1002}),(\ref{eq1003}), and (\ref{eq1015}) we get
\begin{eqnarray}
\label{eq1004} \int_0^T |  z'''(L,t)|^2 + |v_{\mc O}''(L,t)|^2 dt \le 2 \left(c_2 +  \frac{C_3({\bf{ G}}_E)}{\sqrt{\lambda_{k'}}}   \right)  \mc{E}_d(0).
\end{eqnarray}
Now if we use
\begin{eqnarray}
\label{eq1016} | \hat{z}'''|^2  \le 2 | z'''|^2 + 2| z_f'''|^2,\quad   | \hat{v}_{\mc O}''|^2  \le 2 |v_{\mc O}''|^2 + 2|  v_{{\mc O}_f}|^2
\end{eqnarray}
together with (\ref{eq1002}) and (\ref{eq1003}), we obtain
    \begin{eqnarray}
\label{eq1007} \int_0^T \left(| z'''(L,t)|^2 + |v_{\mc O}''(L,t)|^2\right)~ dt \ge  \left(\frac{c_1}{2} - \frac{C_{3}({\bf{ G}}_E)}{2\sqrt{\lambda_{k'}}} \right)  \mc{E}_d(0).
\end{eqnarray}
 Therefore for $T>\tau$ inequalities (\ref{eq1004}) and (\ref{eq1007}) give
     \begin{eqnarray}
 \nonumber \left(\frac{c_1}{2} - \frac{C_{3}({\bf{ G}}_E)}{2\sqrt{\lambda_{k'}}} \right) \mc{E}_d(0)\le \int_0^T \left(| z'''(L,t)|^2 + |v_{\mc O}''(L,t)|^2\right)~ dt \le  2 \left(c_2 +  \frac{C_3({\bf{ G}}_E)}{\sqrt{\lambda_{k'}}}   \right)  \mc{E}_d(0)\quad
\end{eqnarray}
By choosing $k'$ large enough as in the assumption together with using (\ref{energy2}), we obtain
  \begin{eqnarray}
 \nonumber \frac{c_1}{2} \mc E(0)\le \left(\int_0^T |  z'''(L,t)|^2 + |v_{\mc O}''(L,t)|^2\right)~ dt \le 2 c_2 C   \mc{E}(0).
\end{eqnarray}
Hence, condition (i) of Theorem  \ref{Haraux-type-thm} is fulfilled.  Condition (ii)  follows from Lemma \ref{overtrivial}.    \vspace{0.1in}

 In the case of (c-D) boundary conditions, (\ref{important-ineq}) takes of the following form
\begin{eqnarray}
\nonumber \|(z, v_{\mc O})^{\rm T}\|^2_{\mc{X}_1}\ge \lambda_{k'} \|(z, v_{\mc O})^{\rm T}\|^2_{\mc{X}}\end{eqnarray}
which means
\begin{eqnarray}
\nonumber \|(z, v_{\mc O})^{\rm T}\|^2_{H^2_0(\Omega)\times(H^1_0(\Omega))^{(m+1)}}\ge \lambda_{k'} \|(z, v_{\mc O})^{\rm T}\|^2_{H^1_0(\Omega) \times (L^2(\Omega)^{(m+1)})}.
\end{eqnarray}
In the case of (m-m) boundary conditions, we use (\ref{evestimate1}) so that (\ref{important-ineq}) takes of the following form
\begin{eqnarray}
\nonumber \|(z, v_{\mc O})^{\rm T}\|^2_{H^1_0(\Omega)\times L^2(\Omega)}\ge \lambda_{k'} \|(z, v_{\mc O})^{\rm T}\|^2_{(L^2(\Omega)/\mathrm{H})\times ((H^1_\dag(\Omega))')^{(m+1)}}.\end{eqnarray}
The rest of the proof  for (c-D) and (m-m) boundary conditions works the same way modulo the obvious modifications. $\square$ \vspace{0.1in}

\subsection{Proof of main observability result}     \label{sub3}     In this subsection we prove our main observability result Theorem \ref{observability}.
We show that the general damped system  is  a bounded perturbation of the undamped system (with $\tilde {\bf{G}}_E = 0$) and if
$\|\tilde {\bf{G}}_E\|$ is sufficiently small, the observability inequalities (Lemma \ref{aaa}) for the  undamped case remain valid.     \vspace{0.1in}

   We will need the the following lemma.
\vspace{0.1in}
   \begin{Lem}\label{ozkan}
Let $T>0.$ For all  $\|\tilde {\bf{G}}_E\|$  sufficiently small  there exists a  constant $C(\tilde {\bf{G}}_E)>0$ such that for all $ t\in (0,T]$
\begin{subequations}
\label{energy101}
\begin{empheq}[left={\phantomword[l]{}{ } \empheqlbrace}]{align}
\label{energy102}   C(\tilde {\bf{G}}_E) ~\mc{E}(0)\le \mc{E}(T)\le \mc{E}(t)\le \mc{E}(0) && \quad\text{(h-N)}, \text{(c-D)} \\
\label{energy103}   C(\tilde {\bf{G}}_E)  ~\mc{E}_{-1}(0)\le \mc{E}_{-1}(T)\le \mc{E}_{-1}(t)\le \mc{E}_{-1}(0) && \quad \text{(m-m),}
\end{empheq}
\end{subequations}
where $\mc E$ and $\mc E_{-1}$ are defined by (\ref{energy}) and (\ref{ext-energy}), respectively.
\end{Lem}
\vspace{0.1in}

\textbf{Proof:} For the (h-N) case, we multiply the first equation in (\ref{mainhomo}) by $\dot z''$ and the second equation in (\ref{mainhomo}) by $\dot v_{\mc O}'',$ and integrate by parts in space and time. For the (c-D) and (m-m) cases, we multiply the first equation in (\ref{mainhomo}) by $\dot z,$ and the second equation in (\ref{mainhomo}) by $\dot v_{\mc O},$ and integrate by parts in space and time. We obtain the following energy identities
\begin{subequations}
\begin{empheq}[left={\phantomword[l]{}{ } \empheqlbrace}]{align}
\nonumber \mc{E}(T)=  \mc{E}(0) -\int_0^T\left<  \tilde {\bf{G}}_E \dot \phi' , {\bf{h}}_E^{-1}\dot \phi'\right>_{\Omega}~dt &&\quad \text{(h-N)}\\
\nonumber  \mc{E}(T)=  \mc{E}(0) -\int_0^T\left<  \tilde {\bf{G}}_E \dot \phi , {\bf{h}}_E^{-1}\dot \phi\right>_{\Omega}~dt&&\quad \text{(c-D),(m-m)}.
\end{empheq}
\end{subequations}
Since the dissipation term is bounded in the natural energy space, there exists a constant $C_1$ such that
\begin{subequations}
\begin{empheq}[left={\phantomword[l]{}{ } \empheqlbrace}]{align}
 \nonumber \left|-\int_0^T\left<  \tilde {\bf{G}}_E \dot \phi', {\bf{h}}_E^{-1}\dot \phi'\right>_{\Omega}~dt\right|\le  C_1 \|\tilde {\bf{G}}_E\|T \mc{E}(0) && \text{(h-N)}\\
\nonumber \left|-\int_0^T\left<  \tilde {\bf{G}}_E \dot \phi, {\bf{h}}_E^{-1}\dot \phi\right>_{\Omega}~dt\right|\le  C_1  \|\tilde {\bf{G}}_E\|T \mc{E}(0), && \text{(c-D),(m-m).}
\end{empheq}
\end{subequations}
 Therefore,  if   $\|\tilde {\bf{G}}_E\|$  is sufficiently small  so that $C( \tilde {\bf{G}}_E):= 1-C_1 \|\tilde {\bf{G}}_E\| T >0, $ i.e. $\|\tilde {\bf{G}}_E\|<\frac{1}{C_1 T},$  then  for each set of boundary conditions
\begin{eqnarray}C(  \tilde {\bf{G}}_E  )  ~\mc{E}(0)\le \mc{E}(T)\le \mc{E}(t)\le \mc{E}(0). \label{group}\end{eqnarray}
 In particular,   (\ref{energy102}) holds.

   Note that (\ref{group}) implies that if $\|\tilde {\bf{G}}_E\|$ is chosen sufficiently small so that $C(\tilde {\bf{G}}_E) >0,$ the semigroup $\{e^{\mc At}\}_{t\ge 0}$ extends to a C$_0$-group on $\mathbb{R}$ for each set of boundary conditions by Proposition 2.7.4 in \cite{Weiss-Tucsnak}.       This remains true of the semigroup extension defined on $\mc{H}_{-1}$.    In particular,  for the case of  (m-m) boundary conditions,   (\ref{energy103}), and hence also  the characterization of $\mc H_{-1}$ in (3.10) remain valid.
   $\square$
   \vspace{0.1in}

   Now we can prove  our main observability result Theorem \ref{observability}.
\newpage

 \textbf{Proof of Theorem \ref{observability}:}


     Consider the (h-N) case.     We write the solution of (\ref{mainhomo}) in the form
 $$[z,v_{\mc O}]^{\rm T}=[z_f,{v_{\mc O}}_f]^{\rm T} +[\hat{z}, \hat{v}_{\mc O}]^{\rm T},$$
 where $[z_f,{v_{\mc O}}_f]^{\rm T}$ solves
  \begin{equation} \label{out1}\left\{ \begin{array}{l}
m\ddot z-\alpha  \ddot z''+K  z''''-   N^{\rm T} {\bf{h}}_E {\bf{ G}}_E \phi_E' +f(x,t) = 0  ~~~~{\mbox{in}} ~~~~\Omega \times \mathbb{R}^+ ,\\
{\ddot v}_{\mc O} -{\bf{p}}_{\mc O}^{-1} {\bf{E}}_{\mc O}  {v}_{\mc O}''+  {\bf{p}}_{\mc O}^{-1}{\bf{h}}_{\mc O}^{-1}{\bf{B}}^{\rm T} {\bf{ G}}_E \phi_E +f_{\mc O}(x,t) = 0~~ {\rm{on}} ~~ \Omega\times \mathbb{R}^+,
   \end{array} \right.
\end{equation}
 with zero initial data and,
\begin{eqnarray}\label{rhs}[f,f_{\mc O}]^{\rm T}=[-N^{\rm T} {\bf{h}}_E \tilde{\bf{ G}}_E \dot\phi_E',{\bf{p}}_{\mc O}^{-1}{\bf{h}}_{\mc O}^{-1}{\bf{B}}^{\rm T}  \tilde{\bf{ G}}_E \dot\phi_E]^{\rm T},
 \end{eqnarray}
 where  $[\hat{z}, \hat{v}_{\mc O}]^{\rm T}$ solves (\ref{mainhomoud})   with the initial data $ [z^0, v_{\mc O}^0, z^1, v_{\mc O}^1]^{\rm T}.$
 Since   (\ref{rhs}) is a bounded coupling term in $\mc H$,    by  equivalence of energy $\mc E_d \asymp \mc E$ (see Remark \ref{equivalence}),    the
 estimates in  part (a) of Theorem \ref{thmdc} (which apply to the \emph{decoupled} system) remain valid for  (\ref{out1}).    Thus for  any  $T> 0$
 we have
    \begin{eqnarray}
\nonumber &&  \int_0^T  \left(|  z_f'''(L,t)|^2 + |v_{{\mc O}_f}''(L,t)|^2\right)~ dt\\
\nonumber && \quad\quad\le \int_0^T \left(\| N^{\rm T} {\bf{h}}_E\tilde{\bf{ G}}_E {\bf{B}}   \dot v_{\mc O}'\|^2_{L^2(\Omega)} +\| N^{\rm T} {\bf{h}}_E\tilde {\bf{ G}}_E {\bf{h}}_E N   \dot z''\|^2_{L^2(\Omega)} \right. \\
\nonumber &&  \quad\quad\quad \left. + \|{\bf{p}}_{\mc O}^{-1}{\bf{h}}_{\mc O}^{-1}{\bf{B}}^{\rm T} \tilde{\bf{ G}}_E {\bf{h}}^{-1}_E {\bf{B}} \dot v_{\mc O}'  \|^2_{(L^2(\Omega))^{m+1}} + \|{\bf{p}}_{\mc O}^{-1}{\bf{h}}_{\mc O}^{-1}{\bf{B}}^{\rm T} \tilde{\bf{ G}}_E N \dot z''  \|^2_{(L^2(\Omega))^{m+1}}\right)~ dt\quad\\
\label{eq1023}   &&\quad\quad\le  C_4(\tilde {\bf{ G}}_E)  \int_0^T \left(\|\dot z''\|^2_{L^2(\Omega)}+ \| \dot v_{\mc O}'  \|^2_{(L^2(\Omega))^{m+1}}   \right)~dt
\end{eqnarray}
where $C_4(\tilde {\bf{ G}}_E)\to 0$ as $\|\tilde {\bf{ G}}_E\|\to 0.$
\vspace{0.1in}

Next, for $T>\tau$ if we apply  part (b) of Theorem \ref{thmdc}    to  $(\hat { z},\hat{y}_{\mc O})^{\rm T}$.       Hence   there exist $c_1,c_2>0$  for which
 \begin{eqnarray}
\label{eq1021} && c_1\mc{E}(0)\le \int_0^T  \left(|  \hat{z}'''(L,t)|^2 + |\hat{v}_{\mc O}''(L,t)|^2\right)~ dt  \le c_2 \mc{E}(0).
\end{eqnarray}
By using (\ref{eq1015}) together with (\ref{energy2}), (\ref{energy102}), (\ref{eq1023}), (\ref{eq1021}) we get
\begin{eqnarray}
\nonumber \int_0^T \left(|  z'''(L,t)|^2 + |v_{\mc O}''(L,t)|^2\right)~ dt \le 2 \left(c_2 +  C_4(\tilde {\bf{ G}}_E)   \right) \mc{E}(0).
\end{eqnarray}
  Now by using (\ref{eq1016}) together with  (\ref{energy102}), (\ref{eq1023}) and (\ref{eq1021}) we get
    \begin{eqnarray}
\nonumber \int_0^T\left( | z'''(L,t)|^2 + |v_{\mc O}''(L,t)|^2\right)~ dt \ge  \left(\frac{c_1}{2} - C(\tilde {\bf{ G}}_E)C_4(\tilde {\bf{ G}}_E) \right) \mc{E}(0).
\end{eqnarray}
For any fixed $T>\tau,$ the constant $C(\tilde {\bf{ G}}_E)$ is bounded for all sufficiently small $\|\tilde {\bf{G}}_E\|$ (See proof of Lemma \ref{ozkan}). Hence, for sufficiently small $\|\tilde {\bf{G}}_E\|,$  we get the desired observability result (\ref{ohinged}).

The rest of the proof  for (c-D) and (m-m) boundary conditions works the same way modulo the obvious modifications. $\square$

\section{Exact controllability results} \label{C2S4-Control} \vspace{0.3in}
Once continuous observability is established on an appropriate
function space, exact controllability will also hold on an appropriately
defined dual space to the observability space. Here we sketch the procedure for the (h-N) case and indicate the modifications for the (c-D) and (m-m) cases.
\subsection{Proof of Proposition \ref{weaksolution} and Theorem \ref{regularity} for the (h-N) case} We first define the transpositional solution of (\ref{maincont}), (\ref{bdrycont3}) and (\ref{initialcont}).

 By Lemma \ref{SAdjoint},  $\mc{A}^*=-\mc{A}(-\tilde {\bf{ G}}_E)$.
 Hence  the dual backward problem   corresponding to (\ref{maincont}), (\ref{bdrycont3}) and (\ref{initialcont})    is given by
\begin{equation}\left\{ \begin{array}{l}
m\ddot {\hat{z}} -\alpha \ddot {\hat{z}}'' +  K   {\hat{z}}'''' -   N^{\rm T} {\bf{h}}_E \left({\bf{ G}}_E {\hat{\phi}_E} - \tilde{{\bf{ G}}}_E \dot {\hat \phi}_E \right)' = 0 ~~~ {\rm{on}}~~ \Omega\times \mathbb{R}^+ \\
   {\bf{h}}_{\mc O} {\bf{p}}_{\mc O} {\ddot {\hat{v}}}_{\mc O} -{\bf{h}}_{\mc O} {\bf{E}}_{\mc O}  {{\hat{v}}}_{\mc O}'' + {\bf{B}}^{\rm T}  \left({\bf{ G}}_E {\hat{\phi}_E} - \tilde{{\bf{ G}}}_E
    \dot {\hat \phi}_E \right) = 0~~ {\rm{on}} ~~ \Omega\times \mathbb{R}^+ \\
{\mbox{where}}~~   {\bf{B}} { {\hat{v}}}_{\mc O}={\bf{h}}_E {\hat{\phi}_E} -{\bf{h}}_E  N {\hat{z}}'
 \end{array} \right.
\label{maindual}
\end{equation}
 with the  boundary and terminal conditions
\begin{eqnarray}  &&{\hat{z}}(0,t)={\hat{z}}''(0,t)= {\hat{z}}(L,t)=0, {\hat{z}}''(L,t) = 0,\quad {{\hat{v}}}_{\mc O}'(0,t) = {{\hat{v}}}_{\mc O}'(L,t)= 0\label{bdrydual4}\\
&&{\hat{z}}(x,T_1)={\hat{z}}^0(x),  ~~\dot {\hat{z}}(x,T_1)={\hat{z}}^1(x), ~~ {{\hat{v}}}_{\mc O}(x,T_1)= {{\hat{v}}}^0_{\mc O}, ~~ {\dot {\hat{v}}}_{\mc O}(x,T_1)= {{\hat{v}}}^1_{\mc O}.\label{initialdual}
 \end{eqnarray}
 Now we multiply the first and second equations in (\ref{maindual}) by  $w''$ and  $y_{\mc O}''$ respectively where $(w,y_{\mc O})^{\rm T}$ is the solution of non-homogenous equation (\ref{maincont})-(\ref{initialcont}), and then integrate by parts using the boundary conditions (\ref{bdrycont3}) and (\ref{bdrydual4}). Combining these (and using the definitions of ${\psi}_E$ and  ${\hat{\phi}_E}$)  yield
 \begin{eqnarray}
\nonumber && 0=\left[\int_\Omega   \left(\dot {\hat{z}}'' \mc Lw - {\hat{z}}'' \mc L \dot w + {\bf{h}}_{\mc O} {\bf{p}}_{\mc O}   {\dot  {\hat{v}}}_{\mc O}'' \cdot y_{\mc O}  - {\bf{h}}_{\mc O} {\bf{p}}_{\mc O} {\hat{v}}''_{\mc O} \cdot \dot y_{\mc O} +  \tilde{{\bf{ G}}}_E  {\hat \phi}_E' \cdot  {\bf{h}}_E \psi_E'\right) ~ dx\right]_0^{T_1}\\
\label{eq122} &&+ \int _0^{T_1} \left(K {\hat{z}}'''(L,t) M(t) + {\bf{h}}_{\mc O} {\bf{E}}_{\mc O}{{\hat{v}}}_{\mc O}'' (L,t)\cdot   {\bf{g}}_{\mc O}(t)\right) ~dt.
  \end{eqnarray}
 Now let $\hat Y:=({\hat{z}}, {{\hat{v}}}_{\mc O}, \dot{\hat{z}}, \dot {{\hat{v}}}_{\mc O})^{\rm T}$ with  $\hat Y(0)=\hat Y_0=({\hat{z}}^0, {{\hat{v}}}^0_{\mc O},  {\hat{z}}^1, {{\hat{v}}}^1_{\mc O})^{\rm T} \in \mc{H},$ and let
 \begin{eqnarray}\label{firstS}\mc{S}= H^1_0(\Omega)\times (L^2_\perp(\Omega))^{(m+1)}\times L^2(\Omega)\times ((\tilde H^1(\Omega))')^{(m+1)}.
 \end{eqnarray}
 where $L^2_\perp(\Omega)=\{\varphi\in L^2(\Omega): \int_{\Omega} \varphi~dx=0 \}=(\tilde L^2(\Omega))'.$ One can easily prove that the map  $\frac{d^2}{dx^2}: H^2_\perp(\Omega)\to L^2_\perp(\Omega)$  is an isomorphism. Moreover, this extends to isomorphism $\frac{d^2}{dx^2}: H^1_\perp(\Omega)\to  (\tilde H^1(\Omega))'.$     Consequently,
 $ \frac{d^2}{dx^2}:   {\mathcal H} \to {\mathcal S}$ is an isomorphism.

 Define $\mc{F}_{T_1}$ to be the linear functional  on $\mc{H}$ by
\begin{eqnarray}
\nonumber &&\mc{F}_{T_1}(\hat Y_0) = \left<  \left(-\mc L  w^1, - {\bf{h}}_{\mc O} {\bf{p}}_{\mc O}  y_{\mc O}^1, \mc L w^0,  {\bf{h}}_{\mc O} {\bf{p}}_{\mc O}  y_{\mc O}^0\right), \hat Y_0'')\right>_{\mc{S}',\mc{S}} \\
 \nonumber &&  -\int _0^{T_1}\left( K {\hat{z}}'''(L,t) M(t) + {\bf{h}}_{\mc O} {\bf{E}}_{\mc O}{{\hat{v}}}_{\mc O}'' (L,t)\cdot   {\bf{g}}_{\mc O}(t)\right)~ dt \\
  \label{eq124}&& +\left< \left(N^{\rm T}\tilde {\bf{ G}}_E ({\bf{h}}_E N{w^0}'' + {\bf{B}}y_{\mc O}^{0'}), -{\bf{B}}^{\rm T} \tilde {\bf{ G}}_E ( N w^{0'} + {\bf{h}}_E^{-1} {\bf{B}} y_{\mc O}^0), 0,0\right), \hat Y_0''\right>_{\mc{S}',\mc{S}}.
 \end{eqnarray}
 Then (\ref{eq122}) becomes
 \begin{eqnarray}\label{eq123} &&\mc{F}_{T_1}(\hat Y_0)=\left.\left<  \left(-\mc L  \dot w , - {\bf{h}}_{\mc O} {\bf{p}}_{\mc O} \dot y_{\mc O}, \mc L w,  {\bf{h}}_{\mc O} {\bf{p}}_{\mc O}   y_{\mc O}\right), \hat Y''\right>_{\mc{S}',\mc{S}}\right|_{t={T_1}}\\
 \nonumber && + \left.\left< \left(N^{\rm T}\tilde {\bf{ G}}_E ({\bf{h}}_E N{w}'' + {\bf{B}}y_{\mc O}'), -{\bf{B}}^{\rm T} \tilde {\bf{ G}}_E ( N w'+ {\bf{h}}_E^{-1} {\bf{B}} y_{\mc O}),0, 0\right), \hat Y''\right>_{\mc{S}',\mc{S}}\right|_{t={T_1}}. \end{eqnarray}
This identity defines a weak solution of (\ref{maincont})-(\ref{initialcont}); more precisely:
 \vspace{0.1in}
 \begin{Def}\label{def-cont} We say that $(w, y_{\mc O}, \dot w, \dot y_{\mc O})^{\rm T}$ is a solution of  (\ref{maincont})-(\ref{initialcont}) on $[0,T]$ if
  $(w, y_{\mc O}, \dot w, \dot y_{\mc O})^{\rm T}\in C([0,T],\mc{C})$    and (\ref{eq123}) is satisfied for all $T_1\in [0,T]$ and for all
  $\hat Y_0\in \mc{H}$ where $\mc C$ is defined by (\ref{control}). \end{Def}
\vspace{0.1in}

  To see that Def. \ref{def-cont} is fulfilled, first note that by Theorem \ref{observability}, $(\hat z'''(L,\cdot), {{\hat{v}}}_{\mc O}'' (L, \cdot))\in (L^2(0,T))^{(m+2)}.$ Furthermore, since  $\hat Y_0 \in \mc{H},$ by Theorem \ref{eigens},  $\hat Y''(\cdot, T_1)\in \mc{S}$ for all $T_1\in [0,T].$ Therefore, for every $T_1\in [0,T]$ the linear form $\mc{F}_{T_1}$ is continuous on  $\mc{H}.$   Consequently the duality pairing in (\ref{eq123}) uniquely defines the $\left(- \mc L  \dot w ,  - {\bf{h}}_{\mc O} {\bf{p}}_{\mc O} \dot y_{\mc O}, \mc L w, {\bf{h}}_{\mc O} {\bf{p}}_{\mc O}   y_{\mc O}\right)^{\rm T}\in \mc{S}'$  where
   $$\mc S'=H^{-1}(\Omega)\times (\tilde L^2(\Omega))^{(m+1)}\times  L^2(\Omega)\times (\tilde H^1(\Omega))^{(m+1)}.$$     But since $$\mc L :H^2(\Omega)\cap H^1_0(\Omega) \to L^2(\Omega) ~~ {\rm{and}} ~~ \mc L : H^{1}_0(\Omega)\to H^{-1}(\Omega)$$
are isomorphisms it follows that  $\left(w(\cdot,t), y_{\mc O}(\cdot,t), \dot w(\cdot,t), \dot y_{\mc O}(\cdot,t)\right)^{\rm T} \in \mc{C}$ for all $t\in \mathbb{R}.$  One can prove the continuity in time, i.e.,  $\left(w(\cdot,t), y_{\mc O}(\cdot,t), \dot w(\cdot,t), \dot y_{\mc O}(\cdot,t)\right)^{\rm T} \in C([0,T], \mc{C})$ through a standard argument; see  e.g., \cite[Theorem 2.5]{Komornik}. This proves Proposition \ref{weaksolution}.
\vspace{0.1in}

 Now we prove Theorem \ref{regularity} by the HUM method (i.e. see \cite[Chapter 4]{Lions}). To apply HUM we seek the controls of the form $(M(t),{\bf{g}_{\mc O}})=(\hat z'''(L,t),\hat v_{\mc O}''(L,t))$ where $(\hat z, \hat v_{\mc O}) $ is the solution of (\ref{maindual})-(\ref{initialdual}) for $T_1=T.$      By the previous discussion,    the backward problem
\begin{equation}\nonumber \left\{ \begin{array}{l}
m\ddot {{w}} -\alpha \ddot {{w}}'' +  K   {{w}}'''' -   N^{\rm T} {\bf{h}}_E \left({\bf{ G}}_E {{\psi}_E} + \tilde{{\bf{ G}}}_E \dot { \psi}_E \right)' = 0 ~~~ {\rm{on}}~~ \Omega\times \mathbb{R}^+ \\
   {\bf{h}}_{\mc O} {\bf{p}}_{\mc O} {\ddot {{y}}}_{\mc O} -{\bf{h}}_{\mc O} {\bf{E}}_{\mc O}  {{{y}}}_{\mc O}'' + {\bf{B}}^{\rm T}  \left({\bf{ G}}_E {{\psi}_E} + \tilde{{\bf{ G}}}_E       \dot { \psi}_E \right) = 0~~ {\rm{on}} ~~ \Omega\times \mathbb{R}^+ \\
{\mbox{where}}~~   {\bf{B}} { {{y}}}_{\mc O}={\bf{h}}_E {{\psi}_E} -{\bf{h}}_E  N {{w}}'
 \end{array} \right.
\end{equation}
 with boundary and terminal conditions
\begin{eqnarray}\nonumber\left\{ \begin{array}{l}
 {{w}}(0,t)={{w}}''(0,t)= {{w}}(1,t)=0, ~~{{w}}''(L,t) =  \hat z'''(L,t) \\
 {{{y}}}_{\mc O}'(0,t) = 0, ~~  {{{y}}}_{\mc O}'(L,t)=  \hat v_{\mc O}''(L,t)\\
 {{w}}(x,T)=0,  ~~\dot {{w}}(x,T)=0, ~~ {{{y}}}_{\mc O}(x,T)= 0, ~~ {\dot {{y}}}_{\mc O}(x,T)= 0
\end{array} \right.
\end{eqnarray}
 has a unique solution satisfying
  \begin{eqnarray}\nonumber &&\left(-\mc L  \dot {{w}}(\cdot,0), - {\bf{h}}_{\mc O} {\bf{p}}_{\mc O} {\dot {{y}}}_{\mc O}(\cdot,0), \mc L w(\cdot,0),  {\bf{h}}_{\mc O} {\bf{p}}_{\mc O}  y_{\mc O}(\cdot,0)\right)^{\rm T}\\
\nonumber && + \left(N^{\rm T}\tilde {\bf{ G}}_E ({\bf{h}}_E N{ w}''(\cdot,0) + {\bf{B}} y_{\mc O}'(\cdot,0)), -{\bf{B}}^{\rm T} \tilde {\bf{ G}}_E ( N  w'(\cdot,0)+{\bf{h}}_E^{-1} {\bf{B}}  y_{\mc O}(\cdot,0)), 0, 0 \right)^{\rm T}\in \mc{S}'.
  \end{eqnarray}
    Hence, the controllability map $\Lambda:\mc{S} \to \mc{S}'$ defined by
    \begin{eqnarray}\nonumber &&\Lambda(  \hat Y_0'')=\left(-\mc L  \dot {{w}}(\cdot,0) ,  - {\bf{h}}_{\mc O} {\bf{p}}_{\mc O} {\dot {{y}}}_{\mc O}(\cdot,0), \mc L w(\cdot,0), {\bf{h}}_{\mc O} {\bf{p}}_{\mc O}   y_{\mc O}(\cdot,0)\right)^{\rm T}\\
    \nonumber  &&+ \left(N^{\rm T}\tilde {\bf{ G}}_E ({\bf{h}}_E N{ w}''(\cdot,0) + {\bf{B}} y_{\mc O}'(\cdot,0)), -{\bf{B}}^{\rm T} \tilde {\bf{ G}}_E ( N  w'(\cdot,0)+ {\bf{h}}_E^{-1} {\bf{B}}  y_{\mc O}(\cdot,0)), 0, 0\right)^{\rm T}
     \end{eqnarray}
     is continuous from $\mc{S}$ into $\mc{S}'.$ Furthermore, if $ Y_0$ such that
     $$({{w}}(\cdot, 0),  {{{y}}}_{\mc O}(\cdot, 0), \dot{{w}}(\cdot, 0), \dot {{{y}}}_{\mc O}(\cdot, 0))^{\rm T}=(w^0, v_{\mc O}^0, w^1, v_{\mc O}^1)^{\rm T},$$
     then the control $(M(t),{\bf{g}}_{\mc O})=(\hat z'''(L,t),\hat v_{\mc O}''(L,t))$ drives the system (\ref{maincont}) to rest in time $T.$ Therefore, Theorem \ref{regularity} is proved if the surjectivity of the map $\Lambda$ is shown.

Now we choose $(M(t),{\bf{g}_{\mc O}}(t))=(\hat z'''(L,t),\hat v_{\mc O}''(L,t))$ in (\ref{eq124}). Then for $T>\tau $ and for all $\hat Y_0\in \mc{H},$ we have
\begin{eqnarray}
\nonumber  \left<\Lambda (\hat Y_0''), \hat Y_0''\right>_{\mc{S}', \mc{S}}&=& \int _0^T \left(K |{\hat{z}}'''(L,t)|^2 + {\bf{h}}_{\mc O} {\bf{E}}_{\mc O} |{{\hat{v}}}_{\mc O}'' (L,t)|^2\right)~ dt\\
\nonumber &\ge& c_2 \mc{E}(0)\ge c_2 \|\hat Y_0''\|_{\mc{S}}^2
 \end{eqnarray}
 where we used (\ref{ohinged})  with the same constant $c_2.$ Since  $\Lambda$ is a bounded and coercive,  by the
 Lax-Milgram theorem   $\Lambda$ is surjective.    This completes the proof for we complete the proof   of Theorem \ref{regularity} for the (h-N) case.

\subsection{Proofs of Proposition \ref{weaksolution} and Theorem \ref{regularity} for (c-D) and (m-m) cases} \label{Results-other} The proofs for (c-D) and (m-m) cases are similar to the proofs for the (h-N) case with several modifications. For example, we multiply the first equation in (\ref{maindual}) by $w$ and the second equation in (\ref{maindual}) by $y_{\mc O}$ where $(w,y_{\mc O})^{\rm T}$ is the solution of non-homogenous equation (\ref{maincont})-(\ref{initialcont}), and then integrate by parts using the appropriate boundary conditions. Then, the definition of transpositional solution changes as the following
  \begin{eqnarray}\label{eq1123} &&\mc{F}_T(\hat Y_0)=\left.\left<  \left(-\mc L  \dot w , - {\bf{h}}_{\mc O} {\bf{p}}_{\mc O} \dot y_{\mc O}, \mc L w,  {\bf{h}}_{\mc O} {\bf{p}}_{\mc O}   y_{\mc O}\right), \hat Y\right>_{\mc{S}',\mc{S}}\right|_{t=T}\\
 \nonumber && + \left.\left< \left(N^{\rm T}\tilde {\bf{ G}}_E ({\bf{h}}_E N{w}'' + {\bf{B}}y_{\mc O}'), -{\bf{B}}^{\rm T} \tilde {\bf{ G}}_E ( N w'+ {\bf{h}}_E^{-1} {\bf{B}} y_{\mc O}),0, 0\right), \hat Y \right>_{\mc{S}',\mc{S}}\right|_{t=T}. \end{eqnarray}
where the space  $\mc{S}$ is defined as the following
{\small{
 \begin{subequations}\label{semigroupdom3}
 \label{sss}
\begin{empheq}[left={\phantomword[r]{}{\mc{S} = }  \empheqlbrace}]{align}
\label{ScD} & \mc H= H^2_0(\Omega)\times \left(H^1_0(\Omega)\right)^{(m+1)}\times H^1_0(\Omega)\times (L^2(\Omega))^{(m+1)} & \text{(c-D)} \\
\label{Smm}& \mc H_{-1}=H^1_0(\Omega)\times \left(L^2(\Omega)\right)^{(m+1)}\times \left(L^2(\Omega)/\mathrm{H}\right)\times ((H^1_\dag(\Omega))')^{(m+1)} & \text{(m-m).}
\end{empheq}
\end{subequations}}}
 In the above the dual of the space $L^2(\Omega)/\mathrm{H}$ is defined in  Lemma \ref{dualch}.

 Note that (\ref{eq1123}) has  ${\hat Y}$ in the right hand side of the duality pairing whereas  ${\hat Y}''$ appeared in (\ref{eq123}) for the case of (h-N) boundary conditions. However, the duality pairing between $\mc S$ and $\mc S'$ is the same.  This leads to control spaces $\mc C$ defined in (\ref{control1}) and (\ref{control3}) of the same Sobolev order in the cases of (h-N) and (m-m) boundary conditions, as one would expect.

   \vspace{0.1in}

  We indicate below other minor modifications needed for (c-D) and (m-m) cases.
 \vspace{0.1in}

\noindent \textbf{(i)~~ (c-D) case:}    In this case the observability result holds on the  concrete space $\mc H=H^2_0(\Omega)\times \left(H^1_0(\Omega)\right)^{(m+1)}\times H^1_0(\Omega)\times (L^2(\Omega))^{(m+1)}$.   However, as a consequence of the definition of transpositional solution,  the controllability is obtained up to an additive two dimensional space in the velocity component defined in (\ref{spaces}).
To explain this we need the following lemma  which is analogous to Lemmata \ref{lem*}, \ref{dualch}.     Proofs can be found in \cite{Ozer-Hansen} and \cite{O-Hansen1}.
\vspace{0.1in}

\begin{lemma} \label{iso1}   (i) The operator  $\mc L$  is an isomorphism from  $H^2_0(\Omega)$ to $\mathrm{M}^\perp$ where $\mathrm{M}$ is defined by (\ref{spaces}),  (ii)   $(L^2(\Omega)/\mathrm{M})'=\mathrm{M}^\perp,$  where the duality is with respect to the $L^2(\Omega)$ inner product.
 \end{lemma}

By (\ref{ScD})  we have  $\mc S'=H^{-2}(\Omega)\times \left(H^{-1}(\Omega)\right)^{(m+1)}\times H^{-1}(\Omega)\times (L^2(\Omega))^{(m+1)}.$ We see that $\mc L \dot w $ is well-defined at any time as an element of $H^{-2}(\Omega)$ by (\ref{eq1123}). Equivalently, $\left<\dot w, \mc L\psi\right>_{L^2(\Omega)}$ is defined for each $\psi\in H^2_0(0,l)$.  However, the range of $\mc L$ on the restricted space $H^2_0(\Omega)$ is $\mathrm{M}^\perp$ where $\mathrm{M}$ is defined by (\ref{spaces}).   Thus by Lemma  \ref{iso1},    $\dot w$ is well-defined on the quotient space $L^2(\Omega)/\mathrm{M}$.

\vspace{0.1in}

\noindent \textbf{(ii)~~ (m-m) case:} We find a similar phenomenon in (m-m) case but in the reverse sense: the observability result holds on a factor space $\mc H_{-1}=H^1_0(\Omega)\times (L^2(\Omega))^{(m+1)}\times (L^2(\Omega)/\mathrm{H}) \times (H^1_\dag(\Omega)')^{(m+1)},$ while the controllability is obtained on a concrete space defined in (\ref{control}).

By (\ref{Smm}) and Lemma \ref{dualch}, we have $\mc S'=H^{-1}(\Omega)\times \left(L^2(\Omega)\right)^{(m+1)}\times~ \mathrm{H}^\perp \times (H^1_\dag(\Omega))^{(m+1)}.$ Therefore, $\mc L \dot w $ is well-defined since $\mc L: H^1_0(\Omega)\to H^{-1}(\Omega)$ is an isomorphism.  Equivalently, $\dot w \in H^1_0(\Omega)$ for all $T\in \mathbb{R}.$ For the well-posedness of $w$ we investigate the well-posedness of the following term
 \begin{equation}\left<\mc L w (x,T), \hat z(x,T)\right>_{L^2(\Omega)}.\label{dp}\end{equation}
By Lemma \ref{dualch}, when (\ref{dp}) is defined for all $\hat z \in (L^2(\Omega)/\mathrm{H}),$ the term $\mc L w (x,T)$ is uniquely defined in $\mathrm{H}^\perp. $ Therefore, $w$
 is uniquely determined as an element in $H^2_\#(\Omega)$ by Lemma \ref{lem*}. $\square$


\begin{thebibliography}{10}


\bibitem{Hansen3}  S.W. Hansen, {\sl{ Several Related Models for Multilayer Sandwich Plates}}, Mathrmatical Models \& Methods in Applied
Sciences, 14 (2004), pp. 1103--1132.

\bibitem{HO1} S.W. Hansen, O. Imanuvilov, {\sl{Exact controllability of a multilayer Rao-Nakra plate with free    boundary conditions}}, Mathematical Control and Related Fields, 1 (2011), pp. 189-230.

\bibitem{HO2} S.W. Hansen, O. Imanuvilov, {\sl{Exact controllability of a multilayer Rao-Nakra plate with clamped    boundary conditions}}, ESIAM, 17 (2011), pp. 1101-1132.


\bibitem{Rajaram-Hansen3} S.W. Hansen, R. Rajaram, {\sl{Riesz basis property and related results for a Rao-Nakra sandwich beam}}, Discrete and Continuous Dynamical Systems Supplement Vol. (2005), pp. 365--375.


\bibitem{Komornik} V. Komornik, {\sl Exact Controllability and Stabilization: The Multiplier Method }, Wiley, New York, 1994.


\bibitem{Komornik-P} V. Komornik, P. Loreti, {\sl Fourier Series in Control Theory }, Springer-Verlag, New York, 2005.

\bibitem{Lagnese-Lions} J.E. Lagnese, J.-L. Lions, {\sl Modeling Analysis and Control of Thin Plates,} Masson, Paris 1988.

\bibitem{LT1} \newblock I. Lasiecka, R. Triggiani, {\newblock \emph{Exact controllability and uniform stabilization of Kirchhoff plates with boundary controls only in $\left. {\Delta w} \right|_\Sigma,$}} \newblock J. Differential Equations, 93 (1991), pp. 62--101.



\bibitem{LT2} I. Lasiecka, R. Triggiani, {\sl{Factor spaces and implications on Kirchhoff equations with clamped boundary conditions}}, Abstr. Appl. Anal.  6 (8) (2001), pp. 441--488.

\bibitem{LT4}
I. Laisecka, R. Triggiani,   {\sl{Control theory for partial differential equations: Continuous and Approximation Theories, Part 2,}}
Cambridge University Press, Cambridge, 2003.

\bibitem{Lions}
J.L. Lions ~ Exact Controllability, stabilization and perturbations for distributed parameter systems.
 \emph{SIAM Rev.}
 30 (1) (1988), pp. 1--68.

 \bibitem{Mead-Marcus} D.J. Mead, S. Markus {\sl{ The forced vibration of a three-layer, damped sandwich beam with arbitrary boundary conditions,}}
 J. Sound Vibr. 10 (1969), pp. 163--175.




\bibitem{Ozer-Hansen} A. \"{O}zkan \"{O}zer, {\sl{Exact boundary controllability and feedback stabilization for a multilayer Rao-Nakra beam,}}
Ph.D. Thesis, Iowa State University, 2011.
\bibitem{O-Hansen1}  A. \"{O}zkan \"{O}zer, S.W. Hansen, {\sl{Exact controllability of  a Rayleigh beam with a single boundary control}}, Math. Control Signals Syst., 23-1 (2011), pp. 199--222.
\bibitem{Pazy} A. Pazy, {\sl {Semigroups of linear operators and applications to partial differential equations,} } Springer-Verlag, New York, 1983.

\bibitem{Rajaram-Hansen2} R. Rajaram, {\sl{Exact boundary controllability result for a Rao-Nakra sandwich beam}}, Systems Control Lett., 56 (2007), pp. 558--567.


\bibitem{Rao-Nakra} Rao, Y.V.K.S, Nakra, B.C., {\sl{Vibrations of unsymmetrical sandwich beams and plates with viscoelastic cores}},
J. Sound Vibr., 34 (3) (1974), pp. 309-326.

\bibitem{Wang}    J-M Wang, G-Q Xu, S-P Yung, {\sl{Exponential stabilization of laminated beams with structural damping and boundary feedback controls}},  SIAM J. Control Optim., 44 (2005), pp. 1575--1597.

\bibitem{Weiss-Tucsnak} M. Tucsnak, G. Weiss, {\sl{Observation and Control for Operator Semigroups}}, Birkhäuser Verlag, Basel, 2009.

\bibitem {Yan-Dowell} M.J. Yan, E.H. Dowell, {\sl{Governing equations for vibratory constrained-layer damping sandwich plates and beams}}, J. Appl. Mech., 39 (1972), pp. 1041--1046.

\end{thebibliography}
\end{document}